\documentclass[11pt, a4paper]{article}

\usepackage{amsmath,amssymb, amsthm}
\usepackage[utf8]{inputenc}
\usepackage{mathrsfs}
\usepackage{fullpage}
\usepackage{authblk}

\usepackage{graphicx}

\usepackage{wrapfig}

\usepackage{dsfont}
\usepackage{tikz}

\usetikzlibrary{arrows,calc}

\usepackage{verbatim}

\usepackage{hyperref}

\usepackage{xcolor}
\hypersetup{
  colorlinks   = true, 
  urlcolor     = {blue!90!black}, 
  linkcolor    = {blue!90!black}, 
  citecolor   = {red!90!black} 
}

\newcommand{\R}{\mathbb{R}}
\newcommand{\N}{\mathbb{N}}
\newcommand{\Z}{\mathbb{Z}}
\newcommand{\Go}{\Gamma _0}
\newcommand{\0}{\mathbf{0}}
\newcommand{\e}{\mathbf{e}}
\renewcommand{\d}{\mathrm{d}}

\newcommand\blfootnote[1]{
  \begingroup
  \renewcommand\thefootnote{}\footnote{#1}%
  \addtocounter{footnote}{-1}%
  \endgroup
}

\theoremstyle{plain}
\newtheorem{thm}{Theorem}[section]
\newtheorem{prop}[thm]{Proposition}
\newtheorem{lem}[thm]{Lemma}
\newtheorem{cor}[thm]{Corollary}

\theoremstyle{definition}

\newtheorem{con}[thm]{Condition}
\newtheorem{defi}[thm]{Definition}
\newtheorem{rmk}[thm]{Remark}
\newtheorem{conj}[thm]{Conjecture}

\author[1]{Viktor Bezborodov \thanks{Email: \texttt{viktor.bezborodov@univr.it}}} 
\author[1]{
Luca Di Persio \thanks{Email: \texttt{luca.dipersio@univr.it}}}
\author[2]{
 Tyll Kr{\"u}ger \thanks{Email: \texttt{tyll.krueger@pwr.wroc.pl}}}
\author[3]{
Mykola Lebid \thanks{Email: \texttt{mykola.lebid@bsse.ethz.ch}}}
\author[2]{
Tomasz O\.za\'nski \thanks{Email: \texttt{tomasz.ozanski@pwr.edu.pl}}}

\affil[1]{
{The University of Verona}}
\affil[2]{
{Wroclaw University of Technology}}
\affil[3]
{(D-BSSE) ETH Z{\"u}rich}

\title{Asymptotic shape and the speed of propagation
of continuous-time continuous-space birth processes}

\begin{document}

\maketitle

\begin{abstract}
 
 We formulate and prove a shape theorem for a continuous-time 
 continuous-space stochastic growth model
 under certain general conditions. 
 Similarly to the classical lattice growth models
 the proof makes use of the subadditive ergodic theorem.
  A precise expression for the speed of propagation
  is given in the case of a truncated free branching birth rate.

\end{abstract}


\section{Introduction}

\blfootnote{Keywords: \emph{shape theorem,
spatial birth process,
growth model}}

Shape theorems have a long history. Richardson 
\cite{Rich73} proved the shape theorem
for the Eden model. Since then, 
shape theorems have been proven in various settings,
most notably for first passage percolation
and
permanent and non-permanent growth models. 
Garet and Marchand \cite{GM12} 
not only prove a shape theorem for the contact
process in random environment, but also 
have a nice overview of existing results.

Most of literature is devoted 
to discrete-space models. A
continuous-space first passage percolation model 
was analyzed by Howard and Newman \cite{HN97},
see also references therein. A shape theorem
for a continuous-space growth model 
was proven by Deijfen \cite{Dei03},
see also Gou{\'e}r{\'e} and Marchand \cite{GM08}.
Our model is naturally connected to 
that model, see the end of Section 2.

{Questions addressed in this article are motivated 
 not only by probability theory
 but also by studies in natural sciences. In particular, 
one can mention a demand to incorporate spatial information in the
description and analysis of 
1) ecology 2) bacteria populations 3) 
tumor growth 4) epidemiology 5) phylogenetics among others, 
see e.g. \cite{Nature}, \cite{treloar2013multiple}, \cite{vo2015melanoma},
and
\cite{tartarini2015adult}. 
Authors often  emphasize that it is preferable 
to use the continuous-space spaces
$ R ^ 2 $ and 
$ R ^ 3 $ as the basic, or `geographic' space, see e.g. \cite{vo2015melanoma}. 
More on connections between theoretical studies and applications can be 
found in \cite{moller2003statistical}.
}

The paper is organized as follows.
In Section 2 we describe the model
and formulate our results, which 
are proven in Sections 3 and 4. 
Section 5 is devoted to computer simulations
and conjectures.
Technical results, 
in particular on the construction of the process,
are collected in the Section 6.

\section{The model, assumptions and results}

We consider a growth model represented 
by a continuous-time continuous-space
Markov birth process.
Let $\Go$ be the collection of finite subsets of $\R^\d$,
 \[
 \Gamma _0(\R ^\d)=\{ \eta \subset \R ^\d : |\eta| < \infty \},
\]
where $|\eta|$ is the number of elements in $\eta$.
$\Go$ is also called the \emph{configuration space},
or the \emph{space of finite configurations}.

 The evolution of the spatial birth process on $\R^\d$
 admits the following description. 
 Let $\mathscr{B}(X)$ 
 be the Borel $\sigma$-algebra 
 on the Polish space $X$.
 If the system is in state $\eta \in \Go$ 
 at time $t$, then the probability 
 that a new particle appears (a ``birth'') in 
 a bounded set $B\in \mathscr{B}(\R ^\d)$
 over time interval $[t;t+ \Delta t]$ is 
 \[
   \Delta t \int\limits _{ B}b(x, \eta)dx + o(\Delta t),
 \]
 and with probability $1$ no two births happen simultaneously.  
 Here $b: \R ^\d \times \Go \to \R _+$ is some function
 which is called the \emph{birth rate}.
Using a slightly different terminology, we can say that
the rate at which a birth 
 occurs in $B$ is $\int_B b(x, \eta)dx$.
We note that it is conventional to call the function $b$
the
`birth rate', even though it is not a rate in the 
usual sense 
(as in for example `the Poisson  process $(N_t)$ has unit jumps at
rate $1$ meaning that 
$\frac{P\{N_{t+\Delta t} - 
N_{t} = 1\}}{\Delta t} = 1$ as $\Delta t \to 0$')
but rather a version of the Radon--Nikodym derivative
of the rate
with respect to the Lebesgue measure.
 
\begin{rmk}
We characterize the birth mechanism 
by the birth rate $b(x,\eta)$ at each spatial position.
Oftentimes the birth mechanism 
is given
 in terms of contributions 
 of individual particles:
 a particle at $y$,
 $y \in \eta$,
 gives a birth at $x$ 
 at rate $c(x,y, \eta)$
 (often 
 $c(x,y, \eta) =\gamma (y,\eta) k(y, x)$,
 where $\gamma (y,\eta)$
 is the proliferation
  rate
 of the particle at $y$,
 whereas the dispersion kernel $k(y, x)$
 describes the distribution 
 of the offspring),
 see e.g.
Fournier and M{\'e}l{\'e}ard
\cite{FournierMeleard}.
As long as we are not interested 
in the induced genealogical structure,
the two ways of describing
the process are equivalent
under our assumptions.
Indeed, given $c$,
we may set
\begin{equation}
  b(x,\eta) = 
 \sum\limits _{y \in \eta} c(x,y,\eta),
\end{equation}
or, 
conversely,
given $b$, we may set 
\begin{equation}
 c(y,x,\eta) = \frac{g(x-y)}
 {\sum\limits _{y \in \eta} g(x-y)} b(x,\eta),
\end{equation}
where $g: \R ^\d \to (0,\infty)$ is a
continuous
function.
Note that $b$ is uniquely determined 
by $c$, but not vice versa.

 \end{rmk}

 We equip $\Gamma _0$ with the $\sigma$-algebra
 $\mathscr{B}(\Go)$
induced by the sets
\begin{equation} \label{sigmaalgebra}
   \mathbf{Ball}(\eta, r) = \left\{ 
   \zeta \in \Gamma _0 \bigr\vert |\eta| = |\zeta|,
   \text{dist} (\eta,\zeta) < r
   \right\}, \ \ \ \eta \in \Gamma _0, r >0,
\end{equation}
where
$\text{dist}  (\eta,\zeta) = 
\min \left\{ \sum\limits _ {i = 1} ^{|\eta|} |x_i - y_i| \middle| 
\eta = \{x_1, ..., x_{|\eta|} \}, 
\zeta = \{y_1, ..., y_{|\eta|} \}  \right\}$.
For more detail on configuration spaces 
see e.g. 
 R\"ockner and Schied \cite{RockSchied}
or
Kondratiev and Kutovyi \cite{Kondkut}.
In particular, the $\text{dist}$ above
coincides with
 the restriction to the space of finite 
 configurations of the metric $\rho$ used in \cite{RockSchied},
 and the $\sigma$-algebra  $\mathscr{B}(\Go)$ introduced 
 above
 coincides with the $\sigma$-algebra
 from \cite{Kondkut}.

We say that a function $f:\R ^\d \to \R _+$
has an exponential moment if there 
exists $\theta > 0$ such that
\[
 \int _{\R^\d} e^{\theta |x|} f(x) dx < \infty. 
\]
Of course, if $f$ has an exponential moment, then automatically
$f \in L ^1(\R ^\d)$.

\emph{Assumptions on} $b$. We will need several assumptions on the birth rate $b$.
\begin{con}[Sublinear growth]\label{sublinear growth}
 The birth rate $b$ is measurable 
and
 there exists 
a function $a:\R ^\d \to \R _+$
with an exponential moment
such that
\begin{equation} \label{sublinear growth for b} 
 b(x, \eta ) \leq \sum\limits _{y \in \eta } a (x-y).
\end{equation}
\end{con}

 \begin{con}[Monotonicity]\label{mon_con}
  For all $ \eta \subset \zeta$,
  \begin{equation}
 b(x, \eta ) \leq b(x, \zeta), \quad x \in \R ^\d.
\end{equation}

 \end{con}

\noindent The previous condition ensures attractiveness, see below.

\begin{con}[Rotation and translation invariance]
The birth rate
  $b$ is  translation and  rotation invariant:
for every $x, y \in \R ^\d$,  $\eta \in \Go$ and
$M \in \textrm{SO}(\d)$,
\begin{gather*}
  b(x + y, \eta + y ) = b(x, \eta), \\
  b(Mx, M \eta ) = b(x , \eta ).
\end{gather*}
Here $ \textrm{SO}(\d)$ is the orthogonal 
group of linear isometries on $\R ^\d$,
and 
for a Borel set $B \in \mathscr{B}(\R ^\d)$
and $y \in \R ^\d$,
$$
B + y = \{z \mid z = x+y, ~ x \in B\} 
$$
$$
M B = \{z \mid z = Mx, ~ x \in B \}. 
$$

\end{con}

\begin{con}[Non-degeneracy]\label{kook}
 Let there exist 
 $c_0, r >0$ such that
\begin{equation} \label{harlot}
 b(x, \eta) \geq c_0 \quad \textrm{wherever } 
 \min\limits _ {y \in \eta } |x-y| \leq r.
\end{equation}
\end{con}


 \begin{rmk}
  Condition \ref{kook} is used
  to ensure that the system grows at least linearly.
  The condition could be weakened 
  for example as follows: 
  
  \emph{For some  $r_2  > r_1  \geq0$
  and all $x, y \in \R ^\d$, 
    }
  \[
   b(y, \{x\}) \geq c_0  \mathds{1} \{r_1 \leq |x-y| \leq r_2 \}.
  \]
Respectively, the proof would become more intricate.
  
 \end{rmk}

 \begin{rmk}
   If $b$ is like in \eqref{munch} and
   $f$ has polynomial tails, then 
   the result of Durrett \cite{Dur83} suggests that we
   should expect a superlinear growth.
   This is in contrast with Deijfen's model,
   for which Gou{\'e}r{\'e} and Marchand \cite{GM12}
   give a sharp condition
   on the distribution of the outbursts
   for linear or
   superlinear growth.

 \end{rmk}

Examples of a birth rate are

\begin{equation} \label{munch}
b(x,\eta) = \lambda \sum\limits _{y \in \eta } f (|x-y|),
\end{equation}
\noindent and 
\begin{equation}\label{motile}
b(x,\eta) = k \wedge \left( \lambda \sum\limits _{y \in \eta } f (|x-y|) \right),
\end{equation}
\noindent where $\lambda, k$  are positive constants and $f: \R_{+} \to \R_{+}$ is a continuous, non-negative, non-increasing function with compact support.

We denote 
 the underlying probability space by
$(\Omega, \mathscr{F}, P)$.
 Let $\mathscr{A}$ be a sub-$\sigma$-algebra of 
 $\mathscr{F}$.
A random element $A$ in $\Go$
is $\mathscr{A}$-measurable if 
\begin{equation}\label{toad}
 \Omega \ni \omega \to A = A(\omega) \in \Go
\end{equation}
is a measurable map from
the measure space $(\Omega, \mathscr{A})$
to  $(\Go, \mathscr{B}(\Go))$.
Such an 
$A$ will also be called an $\mathscr{A}$-\emph{measurable 
finite
random} set.

The birth process
will be obtained as a unique solution 
to a certain stochastic equation.
The construction
and the proofs of key properties,
such as the rotation invariance and the 
strong Markov property,
are given in Section \ref{sec6}.
We place the construction toward the end 
because it is rather technical
and the methods used there
do not shed much light 
on the ideas of the proofs of our main results.
Denote by $(\eta _t ^{s, A})_{t \geq s}
 = (\eta _t ^{s, A}, t\geq s)$ the process 
started at time $s \geq 0$ from 
an  $\mathscr{S}_s$-measurable finite random set 
$A$.
Here 
 $(\mathscr{S}_s)_{s\geq 0}$ is a filtration
of $\sigma$-algebras to which 
$(\eta _t ^{s, A})_{t \geq s} $ is adapted;
it is introduced after \eqref{nonchalant}. 
Furthermore, $(\eta _t ^{s, A})_{t \geq s} $
is a strong Markov process 
with respect to $(\mathscr{S}_s)_{s\geq 0}$
- see Proposition \ref{strongMP}. 

The construction method we use
has the advantage that the 
stochastic equation approach
resembles graphical representation
(see e.g. Durrett \cite{Dur88}
or Liggett \cite{Liggettbook2})
in the fact that it preserves 
monotonicity:
if $s \geq 0$ and a.s. $A \subset B$,
$A$ and $B$ being 
$\mathscr{S}_s$-measurable finite 
random  sets,
then a.s.
\begin{equation}\label{insolence}
 \eta _{t} ^{s, A} \subset  \eta _{t} ^{s, B}, \ \ \ t \geq s.
\end{equation}
This property
is proven in Lemma \ref{couple}
and is often refered to as 
\emph{attractiveness}.

The process
started from a single particle at $ \0 $ at time zero will be denoted 
by $(\eta _t )_{t \geq 0} $; thus, 
$\eta _t = \eta ^{0, \{ \0 \}} _t$.
Let 
\begin{equation} \label{salvo}
 \xi _t : = \bigcup\limits _{x \in \eta _t} B(x, r)
\end{equation}
and similarly 
\[
 \xi ^{s, A} _t : = \bigcup\limits 
 _{x \in \eta  ^{s, A} _t} B(x, r),
\]
where $B(x, r)$ is the closed ball of radius $r$
centered at $x$ (recall that $r$ appears in \eqref{harlot}).

The following theorem represents the main result of the paper.

\begin{thm} \label{shape thm}
 There exists $\mu >0$ such that
 for all $\varepsilon \in (0,1) $
 a.s. 
  \begin{equation} \label{foolhardy}
   (1- \varepsilon) B(\0, \mu ^{-1}) \subset
   \frac{\xi_t}{t} \subset (1+ \varepsilon) B(\0, \mu ^{-1}) 
  \end{equation}
for sufficiently large $t$.
  
\end{thm}

\begin{rmk}
Let us note that the statement of
Theorem \ref{shape thm} does not depend 
on our choice for the radius in \eqref{salvo}
to be $r$; we could have taken any positive
constant, for example 
\[
 \bigcup\limits _{x \in \eta _t} B(x, 1)
\]
In particular, $\mu$ 
in \eqref{foolhardy} does not depend on $r$.
\end{rmk}

 It is common
to write the ball radius as
the reciprocate $\mu ^{-1}$, probably
because  $\mu $ comes up
in the proof as the limiting value
of a certain sequence of random variables
after applying the subadditive ergodic theorem;
 see 
e.g. Durrett \cite{Dur88} or Deijfen
\cite{Dei03}.
We decided to keep the tradition
not only for historic reasons,
but also because $\mu$
comes up as a certain limit 
in our proof too, even though we do not obtain $\mu$
directly from the subadditive ergodic theorem.
The value $\mu ^{-1}$ is called 
the \emph{speed of propagation}.
The subadditive ergodic 
theorem is a cornerstone 
in the majority of  shape theorem proofs,
and our proof  relies on it.

\emph{Formal connection to Deijfen's model}. The model introduced 
in \cite{Dei03}
with deterministic outburst radius, that is, 
when
in the notation of \cite{Dei03}
the 
distribution of ourbursts $F$ is the Dirac measure: $ F = \delta _R$
for some $R \geq 0$,
can be identified with
\[
 \zeta ^R _t = \bigcup\limits _{x \in \eta _t} B(x, R)
\]
for the birth process $(\eta _t)$
with birth rate 
\[
 b(x,\eta) = \mathds{1} \{\exists y \in \eta : |x-y| \leq R  \}.
\]

\emph{Explicit growth speed for a particular model}.
The precise evaluation of speed appears to be a difficult problem.
For a general one dimensional
branching random walk
the speed of propagation is given by Biggins \cite{Big95}.
An overview of related results 
for different classes of models
can be found in 
Auffinger, Damron, and Hanson
\cite{ADH15}.

Here we give the speed 
for a model with interaction.

\begin{thm}\label{sidekick}
 
Let $\d = 1$ and 

\begin{equation} \label{truncated_b_rate}
 b(x, \eta ) = 2 \wedge \left( \sum\limits _{y \in \eta} \mathds{1}\{ |x-y| \leq 1 \} \right).
\end{equation}

Then the speed of propagation is given by 

\begin{equation}
 \mu ^{-1} = \frac{144\ln (3) - 144 \ln (2) - 40}{25} \approx 0.73548 ...
\end{equation}

\end{thm}

\section{Proof of Theorem \ref{shape thm}} 

We will first show that the 
system grows not faster than linearly.

\begin{prop}\label{incisive}
 There exists $C_{\text{upb}} >0$ such that 
 a.s. for large $t$,
 \begin{equation}\label{forfeiture}
  \eta _t \subset B(\0, C_{\text{upb}} t)
 \end{equation}

\end{prop}
\textbf{Remark}. The index `upb' hints on 
`upper bound'.

\textbf{Proof}. 
It is sufficient to show that 
for $\e = (1,0,...,0) \in \R^\d$
there exists $C>0$ such that
a.s. for large $t$
\begin{equation}\label{licentious}
  \max\{ \langle x, \e \rangle : x\in  \eta _t \} 
  \subset Ct.
 \end{equation}
Indeed, if \eqref{licentious} holds, then by 
Proposition \ref{turd}
it is true if we replace $\e$
with any other unit vector
along any of the $2\d$ directions in $\R^\d$, 
and hence \eqref{forfeiture} holds too.

For $z \in \R$, $y = (y_1,...,y_{\d - 1})
\in \R ^{\d - 1}$ we
define $z \circ y$ 
to be the concatenation $(z, y_1,...,y_{\d - 1})
\in \R ^\d$.
In this proof we denote by $(\bar \eta _t)$
the birth process with
$\bar \eta _0  = \eta _0$ and
the birth rate given 
by the right hand side of \eqref{sublinear growth for b},
namely 
\begin{equation}\label{cognoscenti}
 \bar b(x, \eta ) = 
 \sum\limits _{y \in \eta } a (x-y).
\end{equation}

Since $b(x, \eta ) \leq \bar b(x, \eta )$,
$x\in \R ^\d$, $\eta \in \Go$,
we have by Lemma \ref{couple} a.s.
$\eta _t \subset \bar \eta _t$
for all $t \geq 0$. 
Thus, it is sufficient to prove 
the proposition for $(\bar \eta _t)$.
The  process
$(\bar \eta _t)$ 
with  rate \eqref{cognoscenti}
is in fact a 
continuous-time continuous-space
branching random walk
(for an overview of  branching 
random walks and related topics, see e.g. 
Shi \cite{Shi15}).
Denote by $\bar \eta ^{\e} _t$
the  projection of $\bar \eta  _t$
onto the line
determined by $\e$. 
The process $(\bar \eta ^{\e} _t)$
is itself a branching random walk,
and by 
Corollary 2 in Biggins \cite{Big95},
the position of the rightmost particle 
$X ^\e _t$ of 
$(\bar \eta ^{\e} _t)$ at time $t$ satisfies 
\begin{equation}\label{skittish}
  \lim\limits _{t \to \infty}
 \frac{X ^\e _t}{t} \to \gamma
\end{equation}
for a certain $\gamma \in  (0, \infty)$. 
The conditions from the
Corollary 2 from \cite{Big95}
are satisfied because of Condition \ref{sublinear growth}.
Indeed, 
$(\bar \eta ^{\e} _t)$
is the branching random walk with the birth 
kernel 
\[
 \bar a ^\e (z) = 
 \int\limits _ {y\in\R^{\d - 1} 
 } a(z\circ y)dy,
\]
that is,
$(\bar \eta ^{\e} _t)$ is the a birth process on $\R ^1$
with the birth rate 
\[
 \bar b(x, \eta ) = 
 \sum\limits _{y \in \eta } \bar a ^\e (x-y),
 \ \ \ \ x \in \R, \ \eta \in \Gamma _0(\R).
\]
Note that $a ^\e (z) = a(z)$
if $\d = 1$.
Hence,
in the notation of \cite{Big95}
for $\theta <0$
\[
 m(\theta, \phi) = \int\limits _{\R \times 
 \R_+ } e^{-\theta  z} 
 e^{-\phi \tau} \bar a ^\e (z) dz
 d\tau  = 
 \frac{1}{\phi}
 \int\limits _ {\R } e^{-\theta | z|} \bar a ^\e (z)
dz
=
\frac{1}{\phi}
\int\limits _ {\R } e^{-\theta  |z|} dz 
\int\limits _ {y\in\R^{\d - 1} } a(z\circ y)dy
\]
\[
 =
 \frac{1}{\phi}
\int\limits _ {\R^\d } e^{-\theta  |\langle x, \e \rangle|} 
 a(x)dx \leq 
  \frac{1}{\phi}
\int\limits _ {\R^\d } e^{-\theta  |x|} 
 a(x)dx,
\]
and thus $\alpha(\theta) < \infty$
for a negative $\theta$ satisfying
$\int\limits _ {\R^\d } e^{-\theta  x} a(x)
dx < \infty$ (the functions $m(\theta, \phi)$
and 
$\alpha(\theta)$
are defined in \cite{Big95}
at the beginning of Section 3).

Since \eqref{licentious} follows from \eqref{skittish},
the proof of the proposition is now complete.
\qed


 Next, using a comparison with the Eden model
 (see Eden \cite{Eden}),
 we will show that 
 the system grows not slower than linearly 
 (in the sense of Lemma \ref{play hooky} below).
 The Eden model is 
 a model of tumor growth on the lattice $\Z ^\d$.
 The evolution starts from a single particle at the origin.
 A site once occupied stays occupied forever.
 A vacant site becomes occupied at  rate $\lambda > 0$
 if at least one of its neighbors is occupied.
 Let us mention that this model 
 is closely related to the 
  first passage percolation model,
  see e.g. Kesten \cite{Kes87}
  and
Auffinger, Damron, and Hanson
  \cite{ADH15}. In fact, the two models
  coincide if the passage times (\cite{Kes87})
  have exponential distribution.

 For $z = (z_1, ..., z_\d) \in \Z ^\d$, 
 let $|z|_1 = \sum\limits _{i=1} ^\d |z_i|$. 
 \begin{lem} 
  Consider the Eden model starting from 
  a single particle at the origin.
  Then there exists a constant $\tilde C  >0$
  such that for every $z \in \Z ^\d$ and time 
  $t  \geq \frac{4 e^2}{\lambda ^2 (e -1)^2 } \vee \tilde C |z|_1$,
  \begin{equation}
   P \{ z \text{ is vacant at } t  \} \leq e ^{- \sqrt{t}}. 
  \end{equation}

 \end{lem}
\textbf{Proof}. Let $\sigma _z$
be the time when $z$ becomes occupied. 
Let $v$ be a path on the integer lattice
of length $m = \text{length}(v)$
starting from $\0$ and ending in $z$, so that
$v_0 = \0$,  $v_{m} = z$,
 $v_i \in \Z  ^ \d$
and $|v_i - v_{i-1}| = 1$,
$i = 1,...,m$. 
Define $\sigma(v)$ as the time it takes
for the Eden model to move along the 
path $v$; that is, if 
$v_0,..., v_j$ are occupied,
then a birth can only occur at $v_{j+1}$. 
By construction 
$\sigma(v)$ is distributed as the sum of 
$\text{length}(v)$ independent unit exponentials
(the so called passage times; see 
e.g. \cite{Kes87} or \cite{ADH15}).
We have 
\[
 \sigma _z = \inf\{ \sigma(v): v \text{ 
 is a path from } \0 \text{ to } z\}.
\]
Hence  
$\sigma _z$
is dominated by the sum
of $|z|_1$ independent unit exponentials,
say $\sigma _z \leq Z_1 + ... + Z_{|z|_1}$. 

We have the equality of the events
 \[
   \{ z \text{ is vacant at } t  \} = \{\sigma _z > t \}.
 \]

Note that $E e ^ {\lambda(1 - \frac 1e) Z_1} = e $.
Using Chebyshev's inequality
$P\{ Z > t \} \leq E e ^ {\lambda(1 - \frac 1e) (Z - t)}$,
we get 

\[
  P \{\sigma _z > t \} \leq 
  P \{
  Z_1 + ... + Z_{|z|_1} >t
  \} \leq 
  E \exp \{ \lambda (1 - \frac 1e) (Z_1 + ... + Z_{|z|_1} - t ) \}
\]
\[
 =
 \left[ E e ^ {\lambda (1 - \frac 1e)Z_1} \right]^{|z|_1}
 e ^{- \lambda (1 - \frac 1e)t}
 =
 e^{|z|_1} e ^{- \lambda (1 - \frac 1e)t}.
\]

Since 
\[
 \frac 12 \lambda (1 - \frac 1e) t  \geq \sqrt{t}, 
\]
for $t  \geq \frac{4 e^2}{\lambda ^2 (e -1)^2 }$,
 we may take $\tilde C = \frac {2e}{\lambda (e - 1)}$.
\qed

We now continue to work with the Eden model.
 
 \begin{lem} \label{acquiesce}
  For the Eden model starting from 
  a single particle at the origin,
  there are constants $c_1 , t_0 >0$
   such that
   \begin{equation}
    P \{  \text{there is a vacant site in } B(0,c_1 t) 
    \cap \Z ^\d
    \text{ at } t \} \leq  e^{-\sqrt[4]{t}}, t \geq t_0
   \end{equation}

 \end{lem}

\textbf{Proof}. By the previous lemma for $c_1 < \frac{1}{\tilde C}$,

\[
 P \{  
 \text{there is a vacant site in } B(0,c_1 t) \cap \Z ^\d
 \text{ at } t
 \} 
\]
\[
 \leq \sum\limits _{z \in B(0,c_1 t)\cap \Z ^\d} P \{ z \text{ is vacant at } t  \}
\]
\[
 \leq |B(0,c_1 t)| e^ {-\sqrt{t}},
\]
where
$|B(0,c_1 t)|$ is the number of integer points
(that is, points whose coordinates
are integers)
inside 
$B(0,c_1 t)$. 
It remains to note that $|B(0,c_1 t)|$ grows only polynomially fast
in $t$.
 \qed
  \begin{rmk}\label{nook}
  Let the growth process $(\alpha _t)_{t\geq 0}$ 
  be a $\Z _+ ^{\Z ^\d}$-valued process
  with 
  \begin{equation}
    \alpha (z) \to \alpha (z) + 1 \quad 
    \text{ at rate } \lambda \mathds{1} \big\{\sum\limits _{\substack{y \in \Z ^\d:\\
   |z-y| \leq 1    }}
   \alpha (y) > 0 \big\}, \quad  z\in \Z ^\d, \ \alpha \in \Z _+ ^{\Z ^\d}, 
  \ \sum\limits _{{y \in \Z ^\d
      }}
   \alpha (y) < \infty,
  \end{equation}
 where $\lambda >0$.  
  Clearly,  Lemma \ref{acquiesce}
  also applies to 
 $(\alpha _t)_{t\geq 0}$,
  since it dominates 
  the Eden process.
 \end{rmk}

 Recall that $r$ appears in \eqref{harlot},
 and $(\xi _t)$ is defined in \eqref{salvo}.

\begin{lem} \label{play hooky}
 There are $ c,  s_0 > 0 $ such that 
 \begin{equation}
  P \{ B(\0, cs) \not\subset \xi _s \} 
  \leq  e^ { - \sqrt[4]{s}},
  \quad s \geq s_0.
 \end{equation}
 
\end{lem}

\textbf{Proof}. For $x \in \R ^\d$
let
 $z_x \in \frac{r}{2\d}\Z ^\d$ 
be uniquely determined by
$x \in z_x+ (-\frac{r}{4\d},\frac{r}{4\d}]^\d$.
Recall that $c_0$ appears in Condition \ref{kook}.
Define

\begin{equation}
 \bar b(x,\eta) = c_0
 \mathds{1} \{ z_x \sim z_y 
 \text{ for some } y \in \eta 
 \},
\end{equation}
where $ z_x \sim z_y$ means that $z_x$ and $z_y$
are neighbors on $\frac{r}{2\d}\Z ^\d$.
Let $(\bar \eta _t)_{t \geq 0}$ be the birth process 
with birth rate $\bar b$.
Note that by \eqref{harlot} for every $\eta \in \Go$,
\[
\bar b(x,\eta) \leq b(x,\eta), \
 \ \ x \in \R ^\d,
\]
hence a.s. $\bar \eta _t \subset \eta _t $ by Lemma \ref{couple},
$t \geq 0$.
Then the `projection' process defined by

\[ 
\overline{\overline{\eta}} _t (z) = \sum\limits _{x \in \bar \eta _t} 
 \mathds{1} \{ x \in z + (-\frac{r}{4\d},\frac{r}{4\d}]^\d 
\}, \ \ \ z \in \frac{r}{2\d}\Z ^\d,
\]
is the process $(\alpha _t)_{t\geq 0}$
from
Remark \ref{nook} 
with $\lambda = c_0 \left( \frac{r}{2\d} \right) ^{\d}$
and the `geographic' space  
$\frac{r}{2\d}\Z ^\d$ instead of $\Z ^\d$,
that is, taking values in $\Z _+ ^{\frac{r}{2\d}\Z ^\d}$
instead of $\Z _+ ^{\Z ^\d}$.
Since $\overline{\overline{\eta}} _t (z_x) > 0$
implies that $x \in \xi _t $,
the desired result follows 
from Lemma \ref{acquiesce}
and Remark \ref{nook}.
\qed

 \emph{Notation and conventions}.
 In what follows for $x, y \in \R^\d$
 we define 
 \[
  [x,y]  = \{z \in \R ^\d \mid z = tx + 
  (1-t)y, t \in [0,1]\}.
 \]
 We call $[x,y]$ an interval.
 Similarly, open or half-open intervals 
 are defined, for example 
 \[
  (x,y]  = \{z \in \R ^\d \mid z = tx + 
  (1-t)y, t \in (0,1]\}.
 \]
  We also adopt the convention
$B(x,0) = \{ x \}$.

 For $x \in \R ^\d$
 and $\lambda \in (0,1)$ we define a stopping time $T _\lambda (x)$
 (here and below, all stopping times are considered
 with respect to the filtration 
 $(\mathscr{S} _t)$
 introduced after \eqref{nonchalant})
 by
\begin{equation}\label{gallivant}
 T _\lambda (x) = \inf \{ t>0 : 
 |\eta_t \cap B(x, \lambda |x|)| > 0 
  \}, 
\end{equation}
and
for $x, y \in \R ^\d$, we define 
\begin{equation} \label{walk out on}
 T _\lambda  (x,y) = \inf \left\{ t>T _\lambda(x) :
|\eta_t ^{T _\lambda (x), \{z_\lambda(x)\}} \cap B(y+ z _\lambda(x) - x, \lambda |y-x|)| > 0 
 \right\}- T_\lambda(x),
\end{equation}
where $z_\lambda(x)$
is uniquely defined by
$\{z_\lambda(x)\} = \eta_{_{T_\lambda(x)}}
\cap B(x, \lambda |x|)$. 
Note that $\{z_\lambda(x) \}$
is a $\mathscr{S}_{T _\lambda (x)}$-measurable
finite random set.
Also,
$T_\lambda( \0) = 0$ and
$T _\lambda  (x,x) = 0$ for $x \in \R ^\d$.
To reduce the number of double 
subscripts,
we will sometimes write $z (x)$
instead of $z_\lambda(x)$.

 Since  for $q \geq 1$
$$
\big\{ x_1+x_2: x_1 \in B(x, \lambda |x|),  x_2 \in B((q-1)x, \lambda (q-1) |x|) 
\big\} =  B(qx, \lambda q |x|),
$$
we have by attractiveness (recall \eqref{insolence})
\[
  T_\lambda (q x)   \leq T_\lambda ( x)+ \left( \inf \{ t>0 : 
 |\eta^{T_\lambda ( x), \eta_{T_\lambda ( x)}}_t \cap B(qx, \lambda q |x|)| > 0 
  \} - T_\lambda ( x) \right)
  \]
  \[
  \leq T_\lambda ( x)+ \left( \inf \{ t>0 : 
 |\eta^{T_\lambda ( x), \{z_\lambda (x)\}}_t \cap B( z_\lambda (x) 
  +(q-1)x, \lambda (q-1) |x|)| > 0 
  \} - T_\lambda ( x) \right),
\]
 that is,
  \begin{equation}\label{ire}
 T_\lambda (q x)  
   \leq 
   T_\lambda(x) + T_\lambda(x, qx),
   \quad \quad x \in \R ^\d \setminus \{\0 \}.
  \end{equation}

Note that by the strong Markov property
(Proposition \ref{strongMP}
and Corollary \ref{cor69}),
\begin{equation}\label{stampede}
 T_\lambda(x, qx) \overset{(d)}{=} T_\lambda( (q-1)x).
\end{equation}

The following elementary lemma is used in 
the proof of Lemma \ref{truant}.
\begin{lem}\label{obsequious}
Let  $B_1 = B(x_1, r_1)$ and $B_2 = B(x_2, r_2)$
be two $\d$-dimensional balls.

$(i)$
There exists a constant $c_{\text{ball}}(\d) >0$
depending on $\d$ only such that
 if $B_1$ and $B_2$ are two balls in $\R ^\d$ 
 and  $x_1 \in B_2$
 then 
 \begin{equation}\label{spurious1}
  \text{Vol} (B_1 \cap B_2) \geq c_{\text{ball}}(\d)
 \big( \text{Vol} (B_1 )
  \wedge \text{Vol} ( B_2)\big),
 \end{equation}
where $\text{Vol} (B )$ 
is the $\d$-dimensional volume of 
$B$.
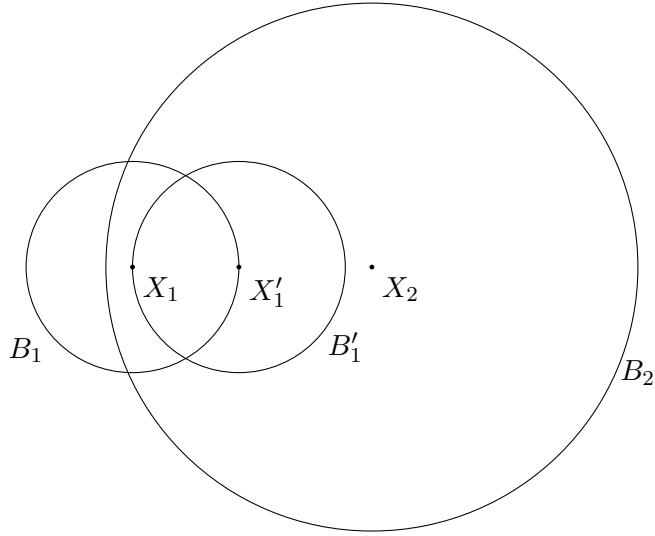
\begin{figure}[!htbp]
	\centering
  \begin{tikzpicture}[scale=0.7]
    \draw (2,2) circle (2.0);
	\filldraw (2,2) node[below right] {$X_1$} circle (1pt);
	\draw (0,0) node [above ] {$B_1$};    
    
    \draw (4,2) circle (2.0);
	\filldraw (4,2) node[below right] {$X'_1$} circle (1pt);    
	\draw (6,0) node [above ] {$B'_1$}; 	

	\draw (6.5,2) node[below right] {$X_2$} circle (5.0);
	\filldraw (6.5,2) circle (1pt);
	\draw (11.5, 0) node {$B_2$};

  \end{tikzpicture}
  \caption{for Lemma \ref{obsequious} $(i)$} \label{fig0}
\end{figure}
$(ii)$ 
The intersection $B_1 \cap B_2$
contains a ball of radius $r _3$ provided that
\[
2 r_3 \leq (r_1 + r_2 - |x_1 - x_2|)\wedge r_1 
\wedge r_2.
\]

\end{lem}
\textbf{Proof}.
$(i)$ Without loss of generality we can
assume that $r_1 \leq r_2$.
Indeed, if $r_1 > r_2$, then 
$x_2 \in B_1$, so we can swap $B_1$ and $B_2$.
Let $B'_1 = B(x'_1, r_1)$ be the shifted ball $B_1$ with 
$x'_1 = x_1 + r_1\frac{x_2 - x_1}{|x_2 - x_1|} $
(see Figure \ref{fig0}). 
The intersection $B'_1 \cap B_1$ is a subset of $B_2$
and is a union of two identical $\d$-dimensional
 hyperspherical caps with height
$\frac{r_1}{2}$. Using the standard formula
for the volume of a hyperspherical cap, we 
see that we can take

\[
c_{\text{ball}}(\d) = \frac{V(B'_1 \cap B_1)}{V(B_1)}=2 \frac{\Gamma(\frac{d}{2}+1)}{\sqrt{\pi} \Gamma(\frac{d+1}{2})}  \int\limits_0^{\frac{\pi}{3}} {\sin^d(s) ds}.
\]
$(ii)$ We have $B_3 \subset B_1 \cap B_2$,
where $B_3 = B(x_3, r_3)$
and $x_3$ is the middle point 
of the interval $[x_1, x_2]
\cap B_1 \cap B_2$.

\begin{lem}\label{truant}
 For every $x \in \R ^\d$ and $\lambda >0$
 there exist $ A_{x, \lambda}, q_{x, \lambda}>0$ 
 such that
   \begin{equation} \label{sublet}
    P \{ T _ \lambda (x) > s \} \leq 
    A_{x, \lambda} e^{-q_{x, \lambda} \sqrt[4]{s}}, \quad s \geq 0.
   \end{equation}

\end{lem}
\textbf{Proof}. 
Let 
\[
 \tau _x = \inf  \{s>0 : x \in \xi_s \}
\]
(recall that $(\xi _t)$ is defined 
in \eqref{salvo}), that is $\tau _x$
is the moment when the first point in the ball $B(x,r)$
appears. 
By Lemma \ref{play hooky} for
$s \geq s_0 \vee \frac{|x|}{c}$
\begin{equation}\label{dither}
 P\{ \tau _x >s \} \leq  P\{ x \notin \xi _s \}
  \leq P\{ B(\0,  |x|) \nsubseteq \xi _s \}
   \leq P\{ B(\0,  cs) \nsubseteq \xi _s \}
 \leq
 e ^{-\sqrt[4]{s}}.
\end{equation}

In the case $r \leq \lambda |x|$
we have a.s. $T _ \lambda (x) \leq \tau _x$,
and
the statement of the lemma follows from 
\eqref{dither}
since for $s \geq s_0 \vee \frac{|x|}{c}$
\[
 P\{ T _ \lambda (x) >s \} 
 \leq P\{ \tau _x >s \}
 \leq
 e ^{-\sqrt[4]{s}}.
\]

Let us now consider the case $r > \lambda |x|$.
Denote by $\bar x \in B(x, r)$ the place
where the particle is born at $\tau _x$.
For $t \geq 0$
on $\{t > \tau _x\}$ we
have 
\[
 \int\limits_{y \in B(x, \lambda |x|)}
 b(y, \eta _t) dy 
 \geq  \int\limits_{y \in B(x, \lambda |x|)}
 b(y, \{\bar x \}) dy 
 \geq 
 \int\limits_{y \in B(x, \lambda |x|)}
 c_0  \mathds{1} \{ y \in B(\bar x, r) \} dy, 
\]
so that by
Lemma \ref{obsequious}
on $\{t > \tau _x\}$
\[
 \int\limits_{y \in B(x, \lambda |x|)}
 b(y, \eta _t) dy 
 \geq 
 \int\limits_{y \in B(x, \lambda |x|)}
 c_0  \mathds{1} \{ y \in B(\bar x, r) \} dy
 \]
 \[
 = c_0  \text{Vol}(B(x, \lambda |x|) \cap 
 B(\bar x, r) )
 \geq c_0 c_{\text{ball}}(\d) \text{Vol}(B(x, \lambda |x|))
 = c_0 c_{\text{ball}}(\d)  V_\d \lambda ^\d |x| ^\d,
\]
where $V_\d =\text{Vol}(B(\0, 1))$,
hence 
\begin{equation*}
\begin{gathered}
 P \{ T _ \lambda (x) - \tau_x > s' \}
 \leq P \{ \inf\{t >0: \eta^{\tau_x, \{\bar x\}}_t
 \cap B(x, r) \ne \varnothing
 \} - \tau_x > s' \}
 \leq e^{ - c_0 c_{\text{ball}}(\d)  V_\d \lambda ^\d |x| ^\d s'}.
 \end{gathered}
\end{equation*}
Combining this with \eqref{dither}
yields the desired result.
\qed

Let us fix an $x \in \R ^\d$, $x \ne \0$, and define 
for $k,n \in \N$, $k < n$,

\begin{equation}
 s_{k,n} = T_\lambda(kx, nx).
\end{equation}

Note that the random variables
$s_{k,n}$ are integrable by Lemma \ref{truant}.
 The conditions of Liggett's subadditive ergodic theorem,
 see \cite{Lig85subadd}, are satisfied here.
 Indeed, condition (1.7) in 
 \cite{Lig85subadd} is 
 ensured by \eqref{ire},
 while conditions (1.8) and (1.9)
 in 
 \cite{Lig85subadd}
 follow from \eqref{stampede} and the strong Markov property
 of $(\eta _t)$ (Proposition \ref{strongMP}
 and Corollary \ref{cor69}).
 Thus, there exists $ \mu _\lambda (x)
 \in [0, \infty)$ such that a.s. 
 and in $L^1$,
 \begin{equation} \label{coax}
  \frac{s_{0,n}}{n} \to \mu _\lambda(x).
 \end{equation}

 \begin{lem} \label{relinquish}
 Let $\lambda >0$.
 For every $x \ne \0$,
 \begin{equation}
  \lim\limits _{t \to \infty} 
  \frac{T _\lambda (tx)}{ t} = \mu _\lambda(x).
 \end{equation}

\end{lem}
\textbf{Proof}.
We know that for every $x \in \R ^\d \setminus \{ \0\}$
 \begin{equation}
  \lim\limits _{n \to \infty} \frac{T _\lambda(nx)}{ n} = \mu
  _\lambda (x).
 \end{equation}
 Denote 
 $\sigma _n = \inf\limits _{y \in [nx, (n+1)x]}T_\lambda(y)$.
 Since there are only a finite number of particles
 born in a bounded time interval, 
 this infinum is achieved. So, 
 let $\tilde z _n$ be such that 
 $\eta_{\sigma _n} \setminus \eta_{ \sigma _n-}
  = \{\tilde z_n \}$. 
  By definition of $\sigma _n$,
  the set 
  \[
   \{y  \in [nx, (n+1)x] \mid 
   \tilde z _n \in B(y, \lambda |y|)\}
  \]
is not empty. $\{ \tilde z _n \}$
is an $\mathscr{S}_{\sigma _n}$-measurable 
finite random set,
so we can apply Corollary \ref{cor69} here.
  
 Define now another stopping time
 \[
  \tilde \sigma _n = \inf\limits 
  \{t>0: \xi^{\sigma _n, \{\tilde z _n\}} _t 
  \supset B(\tilde z _n, \lambda |x| + |x| +2r)\}.
 \]
 Let us show that 
 \begin{equation}\label{hose}
 \sup\limits _{y \in [nx,(n+1)x]} T_\lambda (y) \leq \tilde \sigma _n.
 \end{equation}
 For any  
 $y \in [nx, (n+1)x]$,
 \[
   |y - \tilde z _n| \leq 
    |\tilde z _n - nx| \vee   |\tilde z _n - (n+1)x| 
    \leq \lambda (n+1)|x| + |x|.
  \]
 Therefore the intersection of the balls 
 $B(\tilde z _n, \lambda |x| + |x| +2r)$
 and $B(y, \lambda |y| )$ contains 
 a ball $\tilde B$ of radius $r$ by Lemma \ref{obsequious}, $(ii)$,
 since 
 \[
  \lambda |x| + |x| +2r + \lambda |y| 
  -\lambda (n+1)|x| - |x| \geq \lambda |x|
  + 2r
  +\lambda n|x| - \lambda (n+1)|x| = 2r.
 \]
Since the radius of $\tilde B$ is $r$
and $\xi^{\sigma _n, \{\tilde z _n\}} _{\tilde \sigma _n} \supset 
B(\tilde z _n, \lambda |x| + |x| +2r) 
\supset \tilde B$,
\[
  \eta^{\sigma _n, \{\tilde z _n\}}
  _{\tilde \sigma _n} \cap \tilde B \ne \varnothing,
\]
and hence
\begin{equation}\label{irate}
  \eta
  _{\tilde \sigma _n} \cap \tilde B \ne \varnothing.
\end{equation}
Since $\tilde B \subset B(y, \lambda |y| )$
for all $y \in [n|x|, (n+1)|x|]$,
\eqref{irate} implies \eqref{hose}.
 
 For $q \geq \big(\lambda |x| + |x| +2r\big)\vee
  cs_0$, 
 by Lemma \ref{play hooky}
 \[
  P \{\tilde \sigma _n - \sigma _n  \geq \frac qc \}
  = P \{B(\tilde z _n, \lambda |x| + |x| +2r)
  \nsubseteq
  \xi^{\sigma _n, \{\tilde z _n\}} _{\frac qc  + \sigma _n}
   \} 
   \]
   \[
   \leq 
   P \{B(\tilde z _n, q)
  \nsubseteq
  \xi^{\sigma _n, \{\tilde z _n\}} _{\frac qc  + \sigma _n}
   \} 
   \leq e^ { - \sqrt[4]{\frac qc}},
 \]
 hence 
\begin{equation}
   P \{\tilde \sigma _n - \sigma _n   \geq q' \}
  \leq e^ { - \sqrt[4]{q'}},
  \quad 
  q' \geq \big(\frac{\lambda |x| + |x| +2r}{c}\big)\vee
  s_0.
\end{equation}
By
the Borel--Cantelli lemma 
\[
 P \{\tilde \sigma _n - \sigma _n   > \sqrt{n} 
 \text{ for infinitely many } n\} = 0,
\]
and since $\sigma _n \leq T_\lambda (nx) 
\leq \tilde \sigma _n$, a.s. for large $n$ 
\[
\tilde \sigma _n < T_\lambda (nx)  + \sqrt{n}
\]
and 
\[
\sigma _n \geq T_\lambda (nx)  - \sqrt{n}.
\]
By \eqref{hose}
\[
 \limsup\limits _{n \to \infty} \frac{\sup\limits _{y \in [nx,(n+1)x]} T_\lambda (y)}
 {n} \leq 
 \limsup\limits _{n \to \infty} \frac{\tilde \sigma _n}
 {n} \leq \limsup\limits _{n \to \infty} \frac{T_\lambda (nx) + \sqrt{n}}
 {n} \leq \mu _\lambda (x),
\]
and
\[
 \liminf\limits _{n \to \infty} \frac{\inf\limits _{y \in [nx,(n+1)x]} T_\lambda (y)}
 {n}  =
 \liminf\limits _{n \to \infty} \frac{\sigma _n}
 {n}
 \geq \limsup\limits _{n \to \infty} \frac{T_\lambda (nx) - \sqrt{n}}
 {n} \geq \mu _\lambda (x).
\]
\qed

\begin{lem} The ratio
  $\frac{\mu _\lambda(x)}{|x|}$ 
  in \eqref{coax} does not depend on $x$, $x \ne \0$.
\end{lem}
  \textbf{Proof}. First let us note that
  for every $x \in \R ^\d \setminus \{ \0\}$ and every $q >0$, 
  \begin{equation} \label{brisk}
  \mu _\lambda(x) = \frac{\mu _\lambda (qx)}{q}
 \end{equation}
by Lemma \ref{relinquish}.

On the other hand,
if $|x| = |y|$ then by Proposition \ref{turd}
\begin{equation} \label{clairvoyance}
 \mu _\lambda(x) = \mu _\lambda(y),
\end{equation}
since the distribution of
$(\eta _t)$ is invariant under rotations
and
we can consider $\mu _\lambda(x)$
as a functional 
acting on the trajectory $(\eta_t)_{t \geq 0}$.
The statement
of the lemma 
follows from \eqref{brisk} and \eqref{clairvoyance}.
\qed

Set 
\[
\mu _\lambda: = \frac{\mu _\lambda(x)}{|x|},  \  \ \ 
x \ne \0.
\]
 As $\lambda $
decreases, $T _\lambda (x)$
increases and therefore 
$\mu _\lambda$  increases too. 
Denote 
\begin{equation}
  \mu = 
  \lim\limits _{\lambda \to 0+}
  \mu _\lambda.
\end{equation}

\begin{lem}
 The constants $\mu _\lambda$ and $\mu$
 are strictly positive:
 $\mu _\lambda > 0$,
 $\mu > 0$.
\end{lem}
\textbf{Proof}.
By Proposition \ref{incisive}
for  $x$ with large $|x|$,
\[
 \eta _{\frac{(1-\lambda)|x|}{C_{\text{upb}}}}
 \subset B(\0, (1-\lambda)|x|),
\]
hence for every $\lambda \in (0,1)$ 
for  $x$ with large $|x|$
\[
 T_\lambda (x) \geq  \frac{(1-\lambda)|x|}{C_{\text{upb}}}.
\]
Thus,
\[
 \mu _\lambda \geq \frac{(1-\lambda)}{C_{\text{upb}}}
\]
and
\[
 \mu  = \lim\limits _{\lambda \to 0+}  \mu _\lambda
 \geq \frac{1}{C_{\text{upb}}}.
\]
\qed

\begin{lem}\label{racy}
Let $q, R >0$.
 Suppose that for all $\varepsilon \in (0,1)$  
 a.s. for sufficiently large $n \in \N$
 \begin{equation}\label{semblance}
    \frac{\eta _{qn}}{qn } 
  \subset  (1+ \varepsilon) B(\0, R)
  \ \ \ \ \
  \left(
  (1- \varepsilon)B(\0, R) 
  \subset \frac{\xi _{qn}}{qn } \right).
 \end{equation}
Then for all $\varepsilon \in (0,1)$   
 a.s. for sufficiently large $t \geq 0$
  \begin{equation*}
\frac{\eta _{t}}{t } 
  \subset  (1+ \varepsilon) B(\0, R)
  \ \ \ \ \
  \left(  
 ( 1- \varepsilon)B(\0, R) 
  \subset \frac{\xi _{t}}{t }  \emph{ respectively} \right).
 \end{equation*}
\end{lem}
\textbf{Proof}. We consider the first case only 
-- the proof of the other one is similar.
Since $\varepsilon \in (0,1)$
 is arbitrary, \eqref{semblance}
 implies that 
 for all $\tilde \varepsilon \in (0,1)$ a.s.
 for large $n \in \N$,
 \begin{equation*}
  \frac{\eta _{q(n+2)}}{qn  } \subset (1 + \tilde \varepsilon) B(\0, R).
 \end{equation*}
 
 Since a.s. $(\eta _t)_{t\geq 0}$ is 
 monotonically growing,
it is  sufficient to note that
\begin{equation*}
  \frac{\eta _t}{t} \subset (1 + \varepsilon) B(\0, R) 
  \quad \textrm{if} \quad
  \frac{ \eta _{ \left \lceil{\frac tq}\right \rceil q + q} }
  {\left \lfloor{\frac tq}\right \rfloor q} 
  \subset (1 + \varepsilon) B(\0, R).
 \end{equation*}
\qed

Recall that $c$ is  a constant from Lemma \ref{play hooky}.

\begin{lem} \label{whack}
 Let $\varepsilon  \in (0,1)$. Then
 a.s. 
 \begin{equation}\label{rapport}
  (1 - \varepsilon) B(\0, \mu ^{-1}) 
  \subset \frac{\xi _m}{m }
 \end{equation}
 for large
 $m$ of the form 
$m = (1+\frac{\lambda 
  \mu ^{-1} _\lambda}{c})n$, $n \in \N$.
\end{lem}
 
 \textbf{Proof}.
  Let $\lambda = \lambda _ \varepsilon >0$
 be chosen so small that 
 \begin{equation}
  (1-\varepsilon) \mu ^{-1}
  \leq 
  \frac{ 1-\frac{\varepsilon}{2}}{1 + 
  \frac{\lambda \mu ^{-1} _\lambda}{c} }
  \mu ^{-1} _{\lambda}.
 \end{equation}
 Such a $\lambda$ exists since 
 \[
  \lim\limits _{\lambda \to 0+}
  \frac{  \mu ^{-1} _{\lambda}}{1 + 
  \frac{\lambda \mu ^{-1}  _\lambda}{c} }
 = \mu ^{-1}.
 \]

 Choose a finite sequence of points $\{ x_j, j= 1,...,N \}$
 such that $x _j \in (1 -
 \frac{\varepsilon}{2}) B(\0, \mu _\lambda ^{-1}) $ and
 \[
  \bigcup\limits _{j } 
  B(x_j, \frac{\varepsilon}{4}c) 
  \supset 
  (1 - \frac{\varepsilon}{2}) B(\0, \mu _\lambda ^{-1}).
 \]
 
 Let $\delta > 0$ be so small that
 $(1 + \delta) (1 -  \frac{\varepsilon}{2}) \leq 
 (1 - \frac{\varepsilon}{4}) $. Since a.s.
 \[
  \frac{T _\lambda (n x_j)}{n |x_j|} \to \mu _\lambda,
 \]
 for large $n$ for every $j \in  \{1,...,N \}$
\begin{equation}\label{berate}
  T _\lambda (n x_j)  \leq n |x_j| (1 + \delta) \mu _\lambda \leq 
  n (1 - \frac{\varepsilon}{2}) (1 + \delta) 
  \leq n (1 - \frac{\varepsilon}{4}),
\end{equation}
 so that the system reaches
 the ball $B(n x_j, \lambda n|x_j|)$ before 
 the time
 $n (1 - \frac{\varepsilon}{4})$. 
 Let $Q_n$ be the random event 
 \[\{T _\lambda(nx_j) \leq n (1 - \frac{\varepsilon}{4})
 \text{ for } j = 1,...,N
 \} = 
 \{\eta _{n (1 - \frac{\varepsilon}{4})}
 \cap B(n x_j, \lambda n |x_j|)
 \ne \varnothing, \text{ for } j = 1,...,N \}.\]
 Note that $P(Q_n) \to 1$ by \eqref{berate},
 and even
 \begin{equation}\label{Verona}
  P \{\bigcup\limits _{m \in \N}
\bigcap\limits _{i=m} ^\infty Q_i \}=1.
 \end{equation}
In other words, a.s. for large $i$ all $Q_i$
occur.
 
 Let $\bar z(nx_j)$ be
 defined as
  $ z(nx_j)$ on  $Q_n $
  and as $nx_j$ on the complement 
  $\Omega \setminus Q_n $
  (recall that $ z(x) = z _\lambda(x)$,
  $x \in \R ^\d$, was defined 
  after \eqref{walk out on}).
  The set $ \{ \bar z(nx_j)\}$
  is a finite random 
  \mbox{$\mathscr{S}_{n (1 - \frac{\varepsilon}{4})}$-measurable} set.
 
 Using Lemma \ref{play hooky},
 we will show that  after an additional 
 time interval of length $(\frac{\varepsilon}{4} 
 + \frac{\lambda 
  \mu ^{-1} _\lambda}{c})n$ 
 the entire ball 
 $(1 -  \frac{\varepsilon}{2})n B(\0, \mu _\lambda ^{-1})$
 is covered by $(\xi _t)$, that is,
  a.s. for large $n$
 \begin{equation} \label{myopic}
  (1 -  \frac{\varepsilon}{2})n B(\0, \mu _\lambda ^{-1})
  \subset 
  \xi _{n (1 - \frac{\varepsilon}{4}) 
  + (\frac{\varepsilon}{4} 
 + \frac{\lambda 
  \mu ^{-1} _{_\lambda}}{c})n}
   = 
   \xi _{ n
 + \frac{\lambda n
  \mu ^{-1} _{_\lambda} }{c}}.
 \end{equation}
 Indeed,
 since 
 \[
 B(nx_j, c \frac{\varepsilon}{4} n)
 \subset 
 B(\bar z(nx_j), c \frac{\varepsilon}{4} n + 
 \lambda |x_j|n)
 \subset  B(\bar z(nx_j), c \frac{\varepsilon}{4} n + 
 \lambda \mu ^{-1} _\lambda n),
 \]
 the series
 \[
   \sum\limits _{n \in \N } 
  P 
  \{ B(n x_j , c \frac{\varepsilon}{4} n ) 
  \not\subset
  \xi ^{(n(1-\frac{\varepsilon}{4}),
  \{ \bar z(n x_j)  \})} _{ n
 + \frac{\lambda 
  \mu ^{-1} _{_\lambda} n}{c}} \ \textrm{  for some } j\}
  \]
  \[
  \leq
  \sum\limits _{n \in \N } 
  P 
  \{ B(\bar z(n x_j) , c \frac{\varepsilon}{4} n
  + {\lambda \mu ^{-1} _\lambda}n) 
  \not\subset
  \xi ^{(n(1-\frac{\varepsilon}{4}),
  \{ \bar z(n x_j) \})} _{ n
 + \frac{\lambda 
  \mu ^{-1} _{_\lambda} n}{c}} \ \textrm{  for some } j\}
 \]
converges by Lemma \ref{play hooky}, thus
a.s.
for large $n$,
\begin{equation}\label{redound}
 B(n x_j , c \frac{\varepsilon}{4} n )
  \subset
  \xi ^{(n(1-\frac{\varepsilon}{4}),
  \{ \bar z(n x_j) \})} _{ n
 + \frac{\lambda 
  \mu ^{-1} _{_\lambda} n}{c}},  \ \ \
  j = 1,...,N.
\end{equation}

By
\eqref{Verona} 
a.s. for large $n$

\begin{equation}\label{redound2}
 B(n x_j , c \frac{\varepsilon}{4} n )
  \subset
  \xi ^{(n(1-\frac{\varepsilon}{4}),
  \{  z(n x_j) \})} _{ n
 + \frac{\lambda 
  \mu ^{-1} _{_\lambda} n}{c}},  \ \ \
  j = 1,...,N.
\end{equation}

Hence the choice of $\{x_j, \ j = 1,...,N \}$
and \eqref{redound2} 
yield \eqref{myopic}.
Because  of our choice of $\lambda$,
\[
(1 -  {\varepsilon})n B(\0, \mu  ^{-1})
\subset
\frac{
(1 -  \frac{\varepsilon}{2})}{
(1 + 
  \frac{\lambda \mu ^{-1} _\lambda}{c})}
n B(\0, \mu _\lambda ^{-1}),
\]
which in conjunction with 
\eqref{myopic} implies that
\eqref{rapport} holds a.s. for large
$m$ of the form 
$(1+\frac{\lambda 
  \mu ^{-1} _\lambda}{c})n$,
  where $n \in \N $.
\qed

\begin{lem} \label{implore}
 Let $\varepsilon \in (0,1)$. Then a.s. for large $n \in \N$
 \begin{equation}\label{implore eq}
  \frac{\eta _n}{n } \subset (1 + \varepsilon) B(\0, \mu ^{-1}).
 \end{equation}

\end{lem}

\textbf{Proof}.
Let $\lambda = \lambda _\varepsilon >0$
be so small that 

\begin{equation}\label{rambunctious}
 (1+\frac{\varepsilon}{2}) B(\0, 
 \mu _\lambda ^{-1})
 \subset
  (1+{\varepsilon}) B(\0, 
 \mu  ^{-1})
\end{equation}
 Let $q \in (\varepsilon, \infty)$ and 
 $A$ be the annulus
\begin{equation}\label{rambunctious3}
 A:= (1 +  q) B(\0, \mu _\lambda ^{-1}) \setminus 
  (1 + \frac12 \varepsilon) B(\0, \mu _\lambda ^{-1}),
 \end{equation}
and $\{ x_j, j = 1,...,N \}$
be a finite sequence such that $x_j \in A$ 
and 
 \[
  \bigcup\limits _{j } B(x_j, {\lambda}|x_j|) 
  \supset 
  A.
 \]
 Define $F : = \{ \eta _n \cap n A \ne 
 \varnothing \textrm{ infinitely often}  \}$. 
 On $F$ there exists a (random) $i \in \{1,...,N \}$
 such that the intersection
\begin{equation}
  \eta _n \cap n B(x_{i}, \lambda |x_i|)
\end{equation}
is non-empty infinitely often.
Define also  
\begin{equation}\label{garish}
 F_i :=\left\{ \eta _n \cap n B(x_{i}, \lambda
 |x_i|)
 \ne \varnothing \text{ infinitely often}
 \right\}
\end{equation}
 Note that 
$F\subset \bigcup\limits _{i=1} ^N
F_i$.
 
 On $F_i$ we have 
 \begin{equation*}
  T_{\lambda }(n x_{i}) \leq n
 \end{equation*}
  infinitely often,
hence our choice of $A$ implies
 \[
 \liminf\limits _{n \to \infty} 
 \frac{ T_{\lambda }(n x_i)}{n|x_i|} \leq 
  \liminf\limits _{n \to \infty} \frac{ n}{ (1 +
 \frac 12 \varepsilon) \mu _{\lambda } ^{-1} n}
 = \mu _{\lambda } \frac{ 1 }{ (1 + \frac 12 \varepsilon) }.
\]
 The last inequality and Lemma \ref{relinquish}
 imply that $P(F_i)=0$
 for every $i \in \{1,...,N \}$. Hence $P(F)=0$ too.
 Setting $q  = 2 \mu _\lambda C_{upb} + 1 $,
 so that the radius of the ball on
 the left-hand side of
 \eqref{rambunctious} 
 \[
 q \mu ^{-1} _\lambda > 2C_{upb},
 \]
 by 
 Proposition \ref{incisive}
 and
 the definition of $F$  we get
  a.s. for large $n$,
 \begin{equation}\label{rambunctious2}
   \frac{\eta _n}{n} \subset 
 (1 + \frac12 \varepsilon) B(\0, \mu _\lambda ^{-1}) 
 \end{equation}
 and 
 the statement of the lemma follows from 
 \eqref{rambunctious}
 and 
 \eqref{rambunctious2}.
 \qed
 
 \textbf{Proof of Theorem} \ref{shape thm}.
 The theorem follows from Lemmas \ref{racy},
 \ref{whack}, and \ref{implore}. 
 Note that 
 \begin{equation}\label{implore eq2}
  \frac{\xi _n}{n } \subset (1 + \varepsilon) B(\0, \mu ^{-1}).
 \end{equation}
 is obtained from Lemma \ref{implore}
 by replacing 
 $\varepsilon$ in \eqref{implore eq}
 with
 $\frac \varepsilon 2$.
 
 \qed

 \section{Proof of Theorem \ref{sidekick}}

We precede the proof of Theorem \ref{sidekick}
with
an auxiliary lemma about Markovian functionals
of a general Markov chain. 
 
  Let $(S, \mathscr{B}(S))$ be a Polish (state) space. Consider a
  (time-homogeneous) Markov 
  chain on $(S, \mathscr{B}(S))$ as a family of probability measures
  on $S^ \infty$. Namely,  on the measurable space
  $
  {(\bar \Omega,\mathscr{F}) = (S^ \infty , \mathscr{B}(S ^\infty ))}
 $
  consider a family of probability measures $\{P_s \}_{s \in S}$ 
  such that for the coordinate mappings 
  \begin{align*}
 X_n: \bar \Omega &\rightarrow S, \\
 X_n (s_1,s_2,&...)  = s_n,
\end{align*}
  the process $X := \{X_n \}_{n \in \Z _+ }$ is a Markov chain
  such that for all $s \in S$
  $$
  P_s \{X_0 =s \} =1,
  $$
  $$
  P_s \{ X_{n+m_j}\in A_j, j=1,...,l \mid \mathscr{F} _n \} 
  = P_{X_n} \{ X_{m_j} \in A_j, j=1,...,l  \}.
  $$
  Here $A_j \in \mathscr{B} (S)$,
  $m_j \in \N$, $ l \in \N$, $\mathscr{F} _n = \sigma \{ X_1,...,X_n \}$.
 The space $S$ is separable, hence
 there exists a transition probability 
 kernel
 $Q: S \times \mathscr{B} (S) \rightarrow [0,1]$ such that
 
 $$
 Q(s,A) = P_s \{ X_1 \in A \}, \ \ \ s\in S, \ A \in \mathscr{B} (S).
 $$
  
  Consider a transformation of the chain $X$, $Y_n = f(X_n)$, where
  $f:S\to \R $ is a Borel-measurable function. 
  Here we will give sufficient 
  conditions for 
 $Y = \{Y_n \}_{n \in \Z _+ }$  to be 
  a Markov chain. 
  A very similar question was discussed by
  Burke and Rosenblatt \cite{BRosenblatt}
  for discrete space Markov chains.

 \begin{lem} \label{lumpability}
  Assume that for any bounded Borel function $h: S\rightarrow S$
  
  \begin{equation}\label{insurgent}
   E_s h (X_1) =  E_q h (X_1) \text{\ whenever \ } f(s)= f(q),
  \end{equation}
  Then $Y$ is a Markov chain.
 \end{lem}
 
 \textbf{Remark}. Condition \eqref{insurgent} is the
 equality of distributions of $X_1$
 under two different measures, $P_s$ and $ P_q$.

 \textbf{Proof}. For the natural filtrations of the
 processes $X$ and $Y$ we have an inclusion
 
 \begin{equation} \label{penultimate}
 \mathscr{F} ^X _n \supset \mathscr{F} ^Y _n, \ \ \  n \in \N,
 \end{equation}
since 
$Y$ is a function of $X$. For $k\in \N$
 and bounded Borel functions $h_j: \R \rightarrow \R$,
 $j=1,2,...,k$, 
 
 \begin{equation} \label{prance}
 \begin{split}
     E_s \left[ \prod\limits _{j=1} ^k 
 h_j(Y_{n+j}) \mid \mathscr{F} ^ X _n  \right] = 
 E_{X_n}  &\prod\limits _{j=1} ^k 
 h_j(f(X_{j})) = \\
  \int _S Q(x_0 , dx_1) h_1(f(x_1)) \int _S Q(x_1 , dx_2) h_2(f(x_2))...
 &\int _S Q(x_{n-1} , dx_n) h_n(f(x_n)) \Bigg| _{x_0 = X _n}
 \end{split}
 \end{equation}

To transform the last integral, we introduce a new kernel:
for $y \in f(S)$ chose $x \in S$ with $f(x) = y$,
and then for $B \in \mathscr{B}(\R)$
define
\begin{equation}\label{secretive}
 \overline Q (y, B) = Q (x, f^{-1}(B)).
\end{equation}
The expression on the right-hand side 
 does not depend on the choice of $x$ because of
 \eqref{insurgent}.
 To make the kernel $\overline Q$ defined on 
 $\R \times \mathscr {B} (\R) $, we set 
 
 \[
\overline Q (y, B) = \mathds{1}_{\{0 \in B \}}, \ y \notin f(S).
\]

Then, setting $z_{n} = f(x_{n})$,
we obtain 
from the change of variables formula for the Lebesgue integral
that

$$
\int _S Q(x_{n-1} , dx_n) h_n(f(x_n)) = \
\int _{\R} \overline Q (f(x_{n-1}) , dz_n) h_n(z_n).
$$

Likewise, setting $z_{n-1} = f(x_{n-1})$, we get

\[
\int _S Q(x_{n-2} , dx_{n-1}) h_n(f(x_{n-1}))
\int _S Q(x_{n-1} , dx_n) h_n(f(x_n)) = 
\]
\[ 
 \int _S Q(x_{n-2} , dx_{n-1}) h_n(f(x_{n-1}))
\int _{\R} \overline Q (f(x_{n-1}) , dz_n) h_n(z_n) = 
\]
$$
\int _{\R} \overline Q(f(x_{n-2}) , dz_{n-1}) h_n(z_{n-1})
\int _{\R} \overline Q (z_{n-1} , dz_n) h_n(z_n).
$$
Proceeding further, we obtain
\[
\int _S Q(x_0 , dx_1) h_1(f(x_1)) \int _S Q(x_1 , dx_2) h_2(f(x_2))...
 \int _S Q(x_{n-1} , dx_n) h_n(f(x_n)) =
\]
\[
\int _{\R} \overline Q(z_0 , dz_1) h_1(z_1) \int _{\R} \overline Q(z_1 , dz_2) h_2(z_2)...
 \int _{\R} \overline Q(z_{n-1} , dz_n) h_n(z_n),
\]
where
$z_0 = f (x _0 )$.

Thus,

\[
E_s \left[ \prod\limits _{j=1} ^k 
 h_j(Y_{n+j}) \mid \mathscr{F} ^ X _n  \right] =
 \]
 \[
\int _{\R} \overline Q(f(X_0) , dz_1) h_1(z_1) 
\int _{\R} \overline Q(z_1 , dz_2) h_2(z_2)...
 \int _{\R} \overline Q(z_{n-1} , dz_n) h_n(z_n).
\]

  This equality and  \eqref{penultimate} imply that
  $Y$ is a Markov chain.
  \qed
  
  \begin{rmk} \label{finagle}
  From the proof it follows that 
   $\overline Q$ is the transition probability kernel for 
   the chain $\{f(X_n)\}_{n\in \Z _+}$.
  \end{rmk}

  \begin{rmk}  
   Clearly, this result holds for a Markov chain which is
   not necessarily defined on a canonical state space
   because the property of a process to be a Markov chain
   depends on its distribution only.
  \end{rmk}

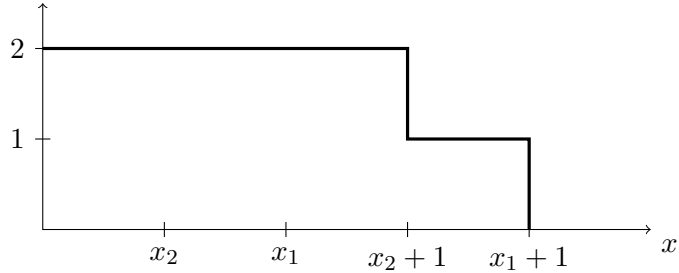
\begin{figure}[!htbp]
\centering
 \begin{tikzpicture}[scale=1.0]
 \coordinate (y) at (0,3);
 \coordinate (x) at (8,0);
 \draw[<->] (y) node[above left] { } -- (0,0) --  (x) node[below right] {$x$};
 
 \coordinate (y1) at ($0.4*(y)$);
 \coordinate (y2) at ($0.8*(y)$);
 
 \coordinate (x1) at ($0.2*(x)$);
 \coordinate (x2) at ($0.4*(x)$);
 \coordinate (x21) at ($0.6*(x)$);
 \coordinate (x11) at ($0.8*(x)$);
 
 \coordinate (vtickshift) at (0.0,0.1);
 \coordinate (htickshift) at (0.1,0.0);
 
 \draw ($(y1)+(htickshift)$) -- ($(y1)-(htickshift)$) node[left] {$1$};
 \draw ($(y2)+(htickshift)$) -- ($(y2)-(htickshift)$) node[left] {$2$};
 \draw ($(x1)+(vtickshift)$) -- ($(x1)-(vtickshift)$) node[below] {$x_2$};
 \draw ($(x2)+(vtickshift)$) -- ($(x2)-(vtickshift)$) node[below] {$x_1$};
 \draw ($(x11)+(vtickshift)$) -- ($(x11)-(vtickshift)$) node[below] {$x_1+1$};
 \draw ($(x21)+(vtickshift)$) -- ($(x21)-(vtickshift)$) node[below] {$x_2+1$};	
 
 \draw [very thick] let \p1=(y2), \p2=(y1), \p3=(x21), \p4=(x11) in
 (\p1) -- (\x3,\y1) -- (\x3,\y2) -- (\x4,\y2) -- (\p4) ; 
 \end{tikzpicture}
  \caption{The plot of $b(\cdot, \eta _t)$.}
\end{figure}
 
 \textbf{Proof of Theorem \ref{sidekick}}.
Without any loss of generality, we will consider the speed of propagation in one direction only,
say toward $+ \infty$. Let $x_1(t)$ and $x_2(t)$ denote the positions 
of the rightmost particle
and the second rightmost particle, respectively ($x_2(t) = 0$ until 
first two births occurs inside $(0, + \infty)$).
Let us observe that
$b(x,\eta _t) \equiv 2$ on $(0,x_2 (t) + 1] $, and
$X = (x_1(t) ,x_2(t))$ is 
a continuous-time pure jump Markov process
on $\{(x_1, x_2)\mid x_1 \geq x_2 \geq 0, x_1 - x_2 \leq 1 \}$ with transition
densities

\begin{equation}
 \begin{aligned}
  (x_1, x_2) &\to (v, x_1) \quad \textrm{ at rate } 1, \quad  v \in (x_2 +1 , x_1 +1 ]; \\
  (x_1, x_2) &\to (v, x_1) \quad \textrm{ at rate } 2, \quad  v \in (x_1 , x_2 +1 ];  \\
  (x_1, x_2) &\to (x_1, v) \quad \textrm{ at rate } 2, \quad  v \in (x_2 , x_1 ].
 \end{aligned}
\end{equation}
(to be precise, the above is true from the moment the first birth inside $\R _+$ occurs).

Furthermore, $z(t) := x_1 (t) - x_2(t)$ 
satisfies 
\[
E \{f(z(t+ \delta)) \mid x_1 (t) = x_1, x_2(t) = x_2  \} = 
E \{f(z(t+ \delta)) \mid x_1 (t) = x_1 + h, x_2(t) = x_2 + h \}
\]
for every $h>0$ and every Borel bounded function $f$. In other words, 
transition rates of $(z(t))_{t \geq 0}$ are entirely determined
by the current state of $(z(t))_{t \geq 0}$,
therefore by Lemma \ref{lumpability} $(z(t))_{t \geq 0}$ is itself
 a pure jump Markov process on $[0,1]$
with the transition densities
\begin{equation}
 \begin{split}
  q(x,y) =  4 & \mathds{1} \{ y \leq x \} + 2 \mathds{1} \{ x \leq y \leq 1-x \} + \mathds{1} \{ y \geq 1-x \}, \quad x \leq \frac12 , y \in [0,1],
  \\
  q(x,y) =  4 & \mathds{1} \{ y \leq 1-x \} + 3 \mathds{1} \{ 1- x \leq y \leq x \} + 
  \mathds{1} \{ y \geq  x \}, \quad x \geq \frac12 , y \in [0,1].
 \end{split}
\end{equation}

Note that the total jump rate out of $x$ is $q(x): = \int _0 ^1 q(x,y) dy = 2+x$.
The process $(z(t))_{t \geq 0}$ is
a regular Harris recurrent Feller process with the Lebesgue measure on $[0,1]$
being a supporting measure
(see e.g. \cite[Chapter 20]{KallenbergFound}). 
Hence the unique invariant measure exists and has a density $g$ with respect
to the Lebesgue measure. 
The equation for $g$ is 
\begin{equation} \label{bland}
 \int\limits _0 ^1 q(x,y) g(x) dx = q(y) g(y).
\end{equation}

Set 
\[
 f (x) = g(x) q(x) \left( \int\limits _0 ^1 g(y)q(y) dy  \right) ^{-1}, \ \ \  x\in [0,1].
\]

It is clear that $f$ is again a density (as an aside we point out that
$f$ is the density of invariant distribution of
the embedded Markov chain of $(z(t))_{t \geq 0}$.
We emphasize however that we do not use this fact in the proof).
Equation \eqref{bland} becomes 
\[
 f(y) = \int\limits _0 ^1 \frac{q(x,y)}{q(x)} f(x)ds,
\]
which after some calculations transforms into

\begin{gather}
 f(y) = 2 \int\limits _0 ^ \frac12 \frac{f(x)dx}{2+x}  \label{croak}
 +
 2 \int\limits _y ^ \frac12 \frac{f(x)dx}{2+x}
 +
 3 \int\limits _\frac12 ^1 \frac{f(x)dx}{2+x}
 +
  \int\limits _\frac12 ^{1-y} \frac{f(x)dx}{2+x}, \quad y \leq \frac 12,
  \\
   f(y) =  \int\limits _0 ^ \frac12 \frac{f(x)dx}{2+x}  \label{fluster}
 +
  \int\limits _0 ^ {1-y} \frac{f(x)dx}{2+x}
 +
  \int\limits _\frac12 ^1 \frac{f(x)dx}{2+x}
 +
 2 \int\limits _y ^1 \frac{f(x)dx}{2+x}, \ \ \quad y \leq \frac 12.
\end{gather}

Differentiating \eqref{croak}, \eqref{fluster}
with respect to
 $y$, we find that $f$ solves the equation

\begin{equation}\label{innuendo}
 \frac{df}{dx} (x) = - 2 \frac{f(x)}{2+x} - \frac{f(1-x)}{3-x}, \quad x \in [0,1].
\end{equation}

Let 
\begin{equation*}
 \varphi (x) := \big[ (2+x)^2(3-x)^2 \big] f(x), \quad x \in [0,1].
\end{equation*}
Then \eqref{innuendo} becomes 

\begin{equation}\label{innuendo2}
 (3-x) \frac{d \varphi}{dx} (x) + 2 \varphi (x) + \varphi (1-x) = 0, 
 \quad x \in [0,1].
\end{equation}

Looking for solutions to \eqref{innuendo2} among polynoms,
we find that $\varphi (x) = c (4 - 3x)$ 
is a solution. By direct substitution
we can check that
\begin{equation}
 f(x) = \frac{c(4-3x)}{(2+x)^2(3-x)^2} \quad x \in [0,1]
\end{equation}
solves \eqref{croak}-\eqref{fluster}.
The constant $c>0$ can be computed, but is irrelevant for our purposes.
Hence, after some more computation,
\begin{equation}\label{vertigo}
 g(x) = \frac{36 (4-3x)}{(2+x)^3(3-x)^2} , \ \ \ x \in [0,1].
\end{equation}

Note that we do not prove analytically that 
equation \eqref{croak}, \eqref{fluster}
has a unique solution. However, uniqueness 
for non-negative integrable solutions
follows from the uniqueness 
of the invariant distribution for $(z(t))_{t \geq 0}$.
Let $l$ be the Lebesgue measure on $\R$.
By an ergodic theorem for Markov processes, see e.g. \cite[Theorem 20.21 (i)]{KallenbergFound},
for any $0 \leq p < p'\leq 1$,
\begin{equation}
 \lim\limits _{t \to \infty} \frac{l \{s: z(s) \in [p,p'], 0\leq s\leq t \}}{t}
 \to \int _p ^{p'} g(x)dx.
 \end{equation}

 Conditioned on $z(t) = z$, the transition densities of $x_1(t)$ are

\begin{equation}
 \begin{aligned}
  x_1 &\to x_ 1 + v \quad \textrm{ at rate } 2, \quad  v \in (0 , 1-z ]; \\
  x_1 &\to x_ 1 + v \quad \textrm{ at rate } 1, \quad  v \in (1-z, 1 ].
 \end{aligned}
\end{equation}

Therefore, the speed of propagation is
\[
 \int\limits _{0} ^1 g(z) dz \left[ \int\limits _{0}^{1-z} 2 y dy + \int\limits _{1-z} ^1 y dy \right]
 =
 \int\limits _{0} ^1 g(z) \left[  (1 - z) ^2  + \frac 12 \big( 1 - (1-z)^2 \big)  \right] dz
\]
\[
 = \int\limits _{0} ^1 g(z)  ( 1 - z +  \frac12 z ^2  ) dz.
\]

Using \eqref{vertigo}, we get

\[
 \int\limits _0 ^1 (1- z +\frac 12 z^2)g(z)dz 
 = \frac{144\ln (3) - 144 \ln (2) - 40}{25}.
\]
\qed

\begin{rmk}
 We see from the proof that the speed can be computed 
 in a similar way for the birth rates of the form
\begin{equation} \label{truncated_b_rate_k}
 b_k(x, \eta) = k \wedge \left( \sum\limits _{y \in \eta} \mathds{1}\{ |x-y| \leq 1 \} \right),
\end{equation}
where 
$k \in  (1,2)$. However, the computations quickly become unwieldy.
 
\end{rmk}

\section{Conjectures and simulations}
\begin{figure}[!htbp]
  \centering
  \begin{minipage}[t]{0.47\textwidth}
    \includegraphics[width=\textwidth]
    {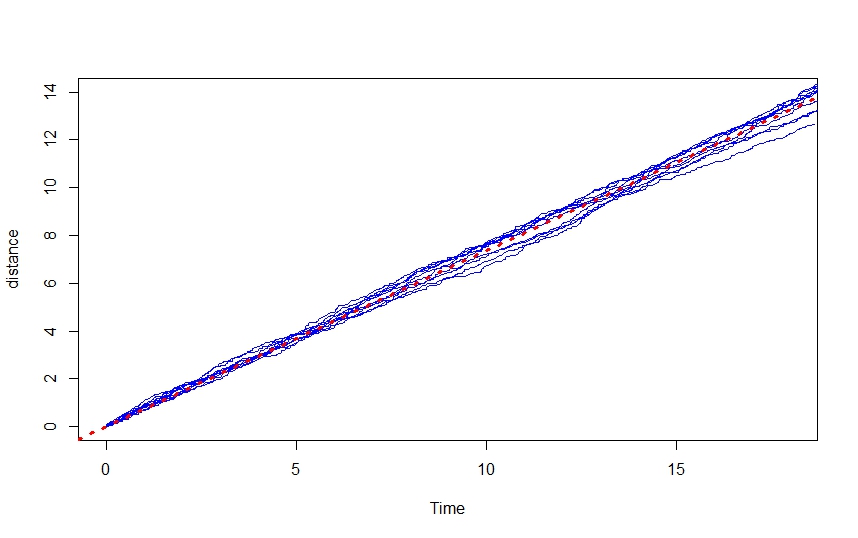}
    \caption{Distance to the furthest particle against time in 
    10 different realizations. Red dashed line shows
    the theoretical distance given by Theorem \ref{sidekick}.} \label{fig1}
  \end{minipage}
  \hfill
  \begin{minipage}[t]{0.47\textwidth}
    \includegraphics[width=\textwidth]{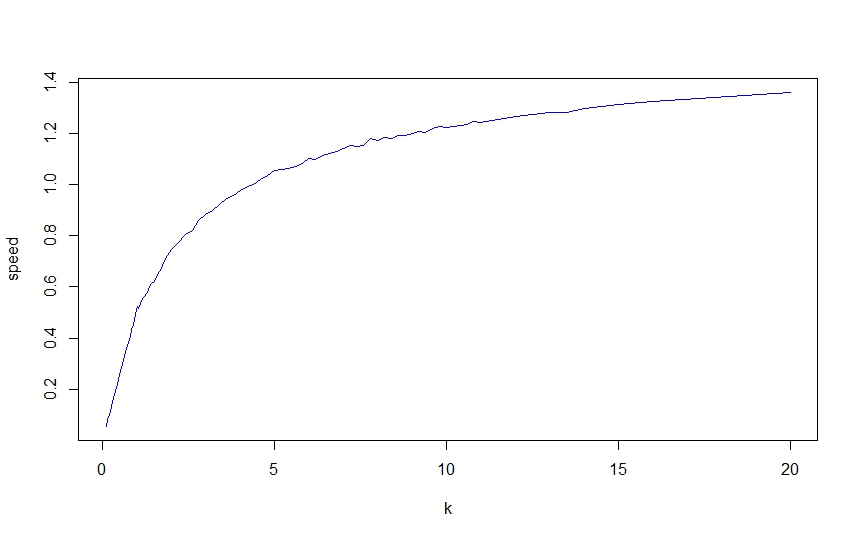}
    \caption{The distance to the furthest particle divided 
    by time against $k$ for the birth rate \eqref{truncated_b_rate_k}
    at time 1000.} \label{fig2}
  \end{minipage}
\end{figure}

In this section we 
collect some conjectures
about the models treated 
in this paper
and related models. 
We also provide 
computer simulations.
{For the latter we mostly use an extended verison 
of the algorithm 
from Section 11.1.2 in \cite{moller2003statistical}}.
We would like to stress that 
this section does not contain
rigorously proven results.
 The models we consider here 
 are mostly one- and two-dimensional.
 The section is divided in five blocks,
 each dealing with its own class of models.

 {\textbf{1}}. Figure \ref{fig1} shows the distance to the 
furthest particle against time for ten different realisations with
birth rate \eqref{truncated_b_rate_k}. 
We see that the fluctuations from 
the value in Theorem \ref{sidekick} 
are rather small. 
Let the dimension $\d = 1$
for the moment.
Denote by $s(k)$
the speed of propagation 
of the system with the birth rate 
\eqref{truncated_b_rate_k}.
Thus, $s(k)$
is $\mu ^{-1}$
in notation of Theorem \ref{shape thm},
if the birth rate is as in \eqref{truncated_b_rate_k}.

Figure \ref{fig2} shows 
 $s(k)$
against $k$ on $x$-axis
for the truncated birth rate
\eqref{truncated_b_rate_k}.

\begin{figure}[!htbp]
  \centering
   \begin{minipage}[t]{0.47\textwidth}
    \includegraphics[width=\textwidth]
    {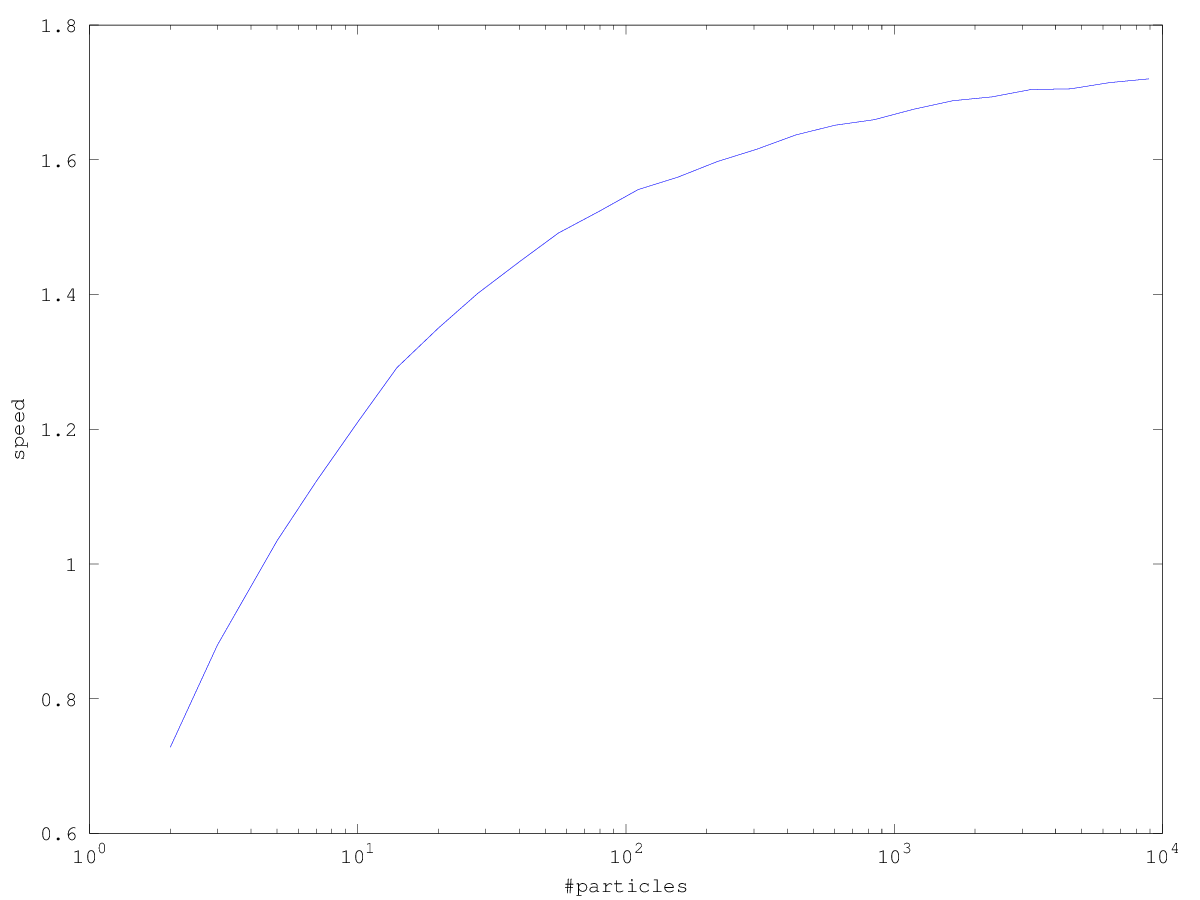}
    \caption{The speed as the distance to the furtherst particle divided 
    by time against the maximal allowed number of particles 
    in the branching random walk with restriction at time $10^4$.
    The scale along $x$-axis is logarithmic. 
    } \label{fig3}
  \end{minipage}
  \hfill
  \begin{minipage}[t]{0.47\textwidth}
    \includegraphics[width=\textwidth]
    {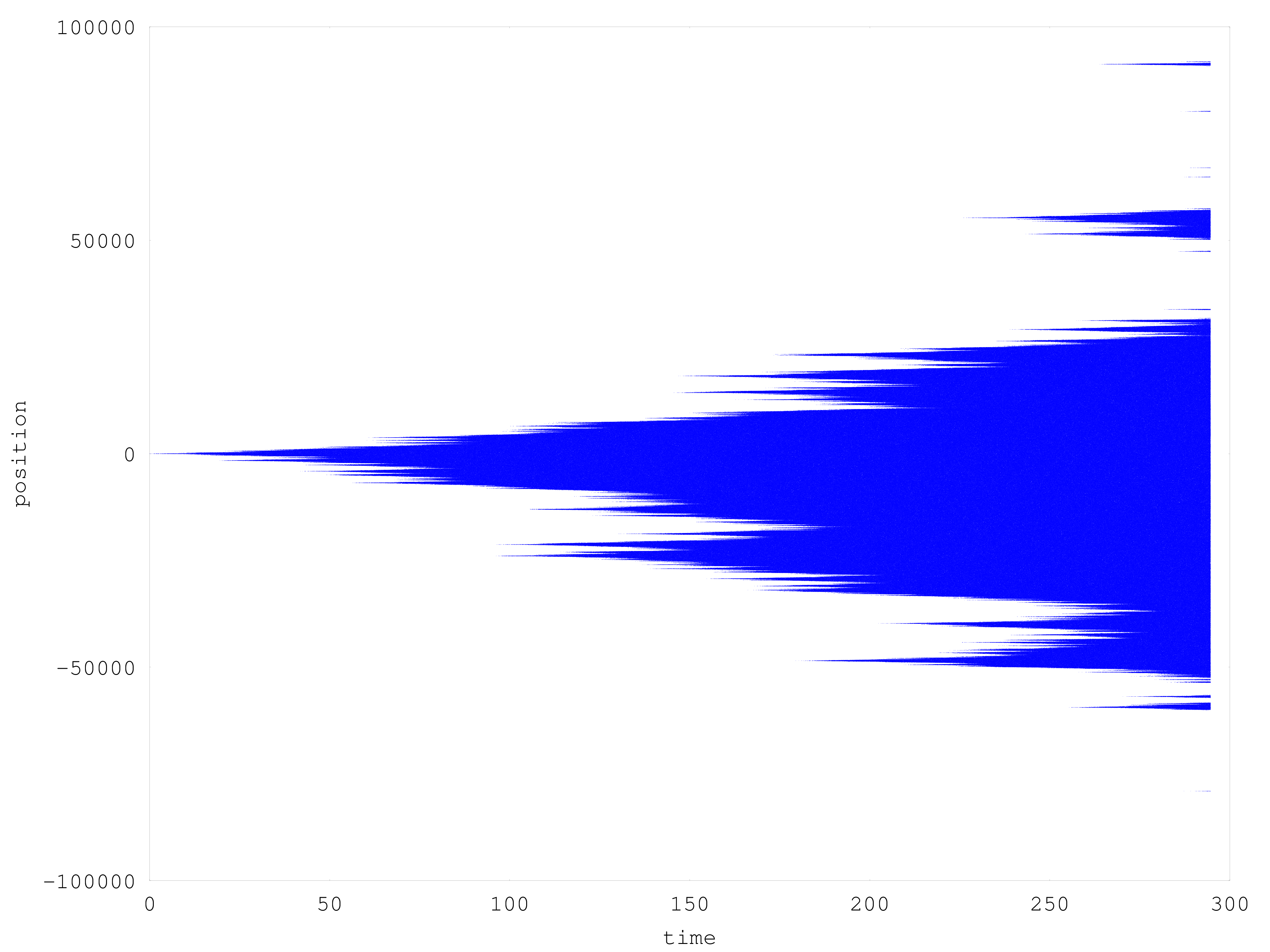}
    \caption{Positions of the occupied sites 
    varying with time for the discrete-space model  
   with birth rate \eqref{aclove} and $\alpha = 2.8$.}
   \label{fig4}
  \end{minipage}
\end{figure}

\begin{figure}[!htbp]
  \centering
  \begin{minipage}[t]{0.47\textwidth}
    \includegraphics[width=\textwidth]
    {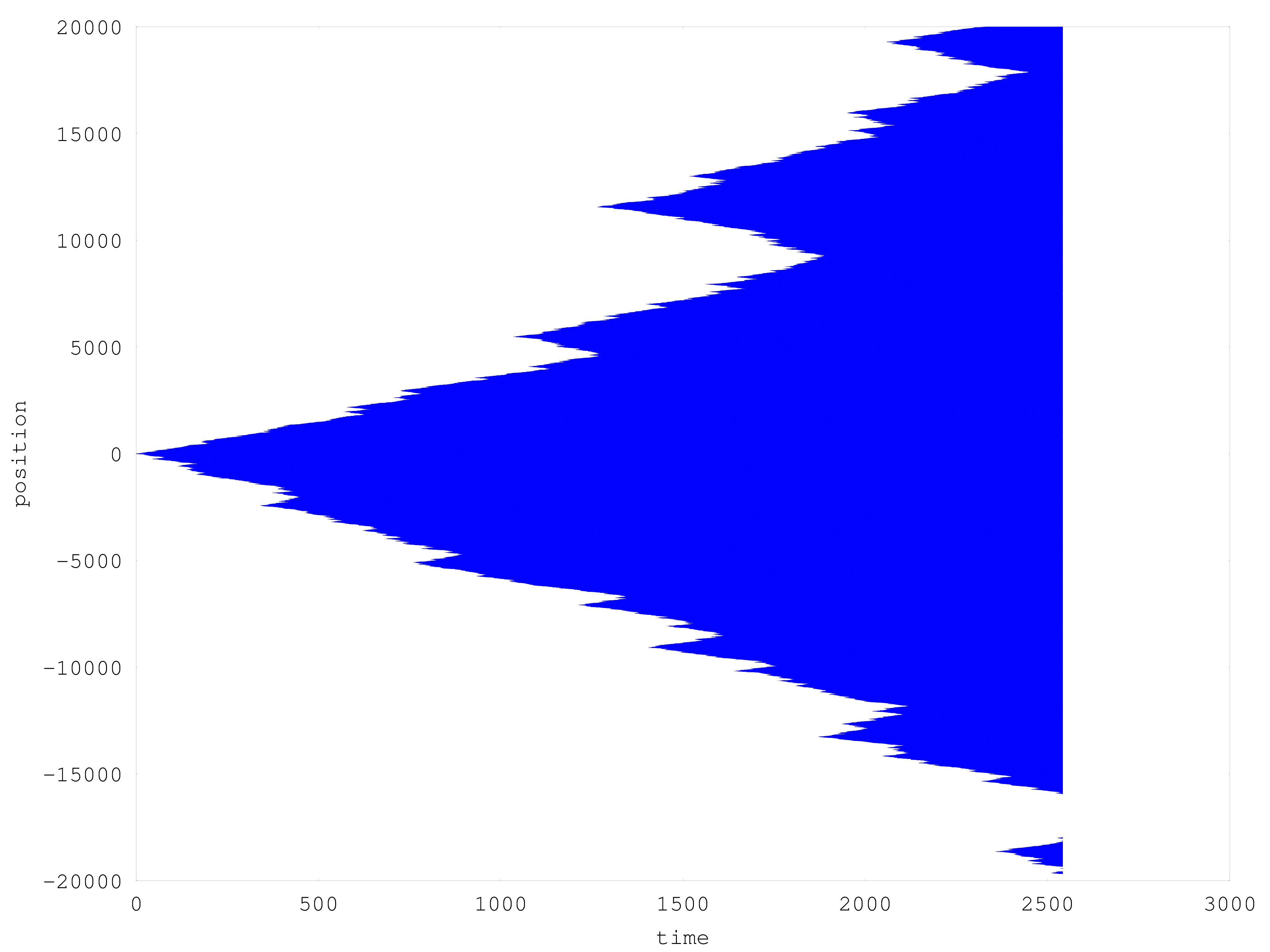}
     \caption{  $\alpha = 3.5$   } \label{fig5}
  \end{minipage}
  \hfill
 \begin{minipage}[t]{0.47\textwidth}
    \includegraphics[width=\textwidth]
    {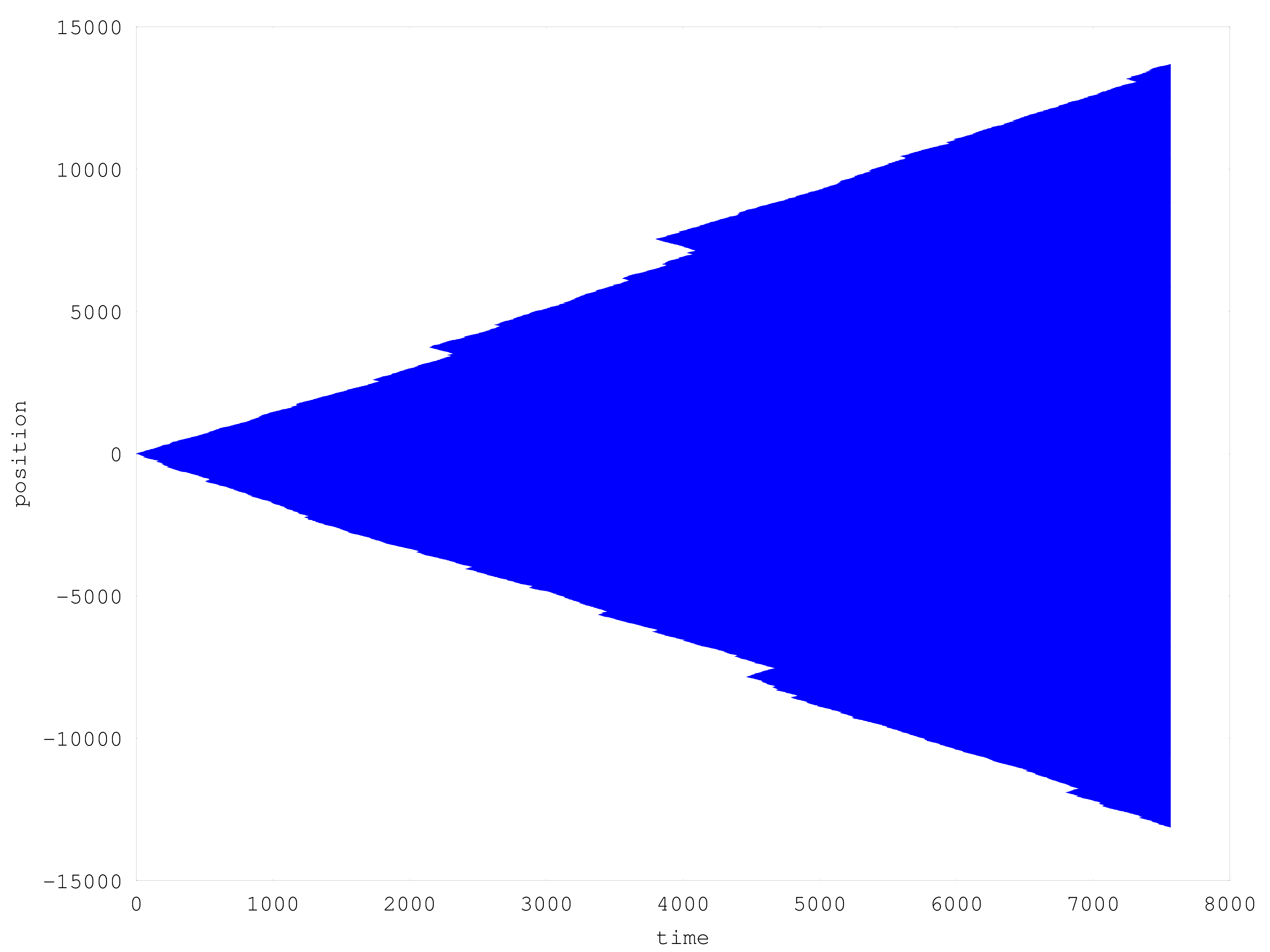}
     \caption{ $\alpha = 4.2$ }
     \label{fig6}
  \end{minipage}
\end{figure}

\begin{conj} 
We conjecture that 
\begin{equation} \label{spurious}
 s(k) \to s _*, \ \ \ k \to \infty.
 \end{equation}
 where $s _*$ is the speed of propagation of the 
 process with the birth rate $b_\infty$ given by
 \begin{equation} \label{free branching indicator}
  b_\infty(x,\eta) =\sum\limits _{y \in \eta} \mathds{1}\{ |x-y| \leq 1 \}.
 \end{equation}
 
 Using the exact formula for 
 for the speed of propagation
of a general branching random walk, see Proposition 1 in Biggins \cite{Big95},
we get 
 \begin{equation}\label{flick}
  s _* = \inf\limits _{a>0} \left\{ 
   \inf\limits _{\theta >0} (e ^{\theta } - e ^{ - \theta } - a \theta ^2) <0
  \right\} \approx 1.81... 
 \end{equation}
The question about the speed of convergence in \eqref{spurious} is more subtle. 

\end{conj}

Clearly, for all $x \in \R ^1$
and $\eta \in \Go$, 
\[
 b_k(x,\eta) \to b_\infty(x,\eta), \ \ \ 
 k \to \infty.
\]
It is probably reasonable to say
that in some sense the birth process
$(\eta _t ^{(\infty)})$
with rate $b_\infty$  from \eqref{free branching indicator}
is approximated by 
the family of processes  $(\eta _t ^{(k)})$
with rate $b_k$ from
\eqref{truncated_b_rate_k}.
Indeed, it can be shown that 
for the stopping time $\theta _k = \inf\{t>0: \eta _t ^{(\infty)}
\ne \eta _t ^{(k)} \}$, we have a.s.
\[
 \theta _k \to  \infty, \ \ \  k\to \infty.
\]
For computational purposes 
$(\eta _t ^{(k)})$
is more tractable since
the number of particles 
within a ball
grows linearly with time,
as opposed to exponential growth for 
$(\eta _t ^{(\infty)})$.

 {\textbf{2}}.
Another way to approximate $(\eta _t ^{(\infty)})$
is to impose restriction on the number of particles.
Specifically, assume that 
the system evolves as the birth process
with rate \eqref{free branching indicator},
but whenever the number of particles 
exceeds a given $n \in \N$,
the leftmost particle is removed.
The number of particles thus stays constant
once it reaches $n$.
Figure \ref{fig3} shows the speed of propagation 
against $n$. The largest $n$ we took was $8902$,
for which the recorded speed is $1.72018$.

 {\textbf{3}}.
Figures \ref{fig4}-\ref{fig6} show the evolution 
of the discrete version of the truncated model
\eqref{motile}: the process evolves 
in 
$\Z ^{\Z _+}$ and the birth rate is 
\begin{equation}\label{aclove}
b(x,\eta) = k \wedge 
\left( \sum\limits _{y \in \eta } a_{\text{pow}} (x-y) \right)
\end{equation}
with 
\[
 a_{\text{pow}} (x) =c_{\text{pow}} \frac{1}{(|x|+1)^{\alpha} }, 
 \ \ \ x \in \Z \setminus \{ 0 \},
\]
\[
 a_{\text{pow}} (0) = 2 c_{\text{pow}},
\]
where $\alpha >2$ and $c_{\text{pow}} 
= c_{\text{pow}}(\alpha)$ is
the normalizing constant.
We have $\alpha = 2.8$ on Figure \ref{fig4}, $\alpha = 3.5$
on Figure \ref{fig5},
and 
$\alpha = 4.2$ on Figure \ref{fig6}.
These pictures allow us to 
observe the development
of the set of occupied sites.
We see that even for a large time,
the set of occupied sites is not a connected interval for $\alpha = 2.8$,
whereas the picture appears to be rather smooth
for $\alpha = 4.2$. We conjecture that
the speed of propagation is superlinear for 
$\alpha = 2.8$, but is linear for $\alpha = 4.2$.
We intend to give a proof in a forthcoming paper.

\begin{figure}[!htbp]
  \centering
  \begin{minipage}[t]{0.47\textwidth}
       \includegraphics[width=\textwidth]
    {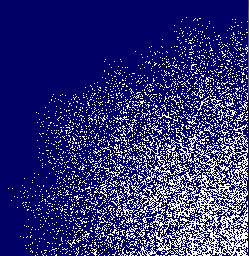}
     \caption{The two dimensional
    simulation of the birth process with rate 
    \eqref{truncated_b_rate_k}, $k = 5$.
    The number of particles is $65000$.} \label{fig7}
  \end{minipage}
  \hfill
 \begin{minipage}[t]{0.47\textwidth}
    \includegraphics[width=\textwidth]
    {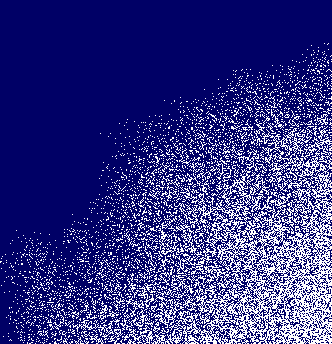}
     \caption{The same system 
     at a later moment.
    The number of particles is $3 \cdot 10^6$.}\label{fig8}
  \end{minipage}
\end{figure}

 {\textbf{4}}.
Figures \ref{fig7} and \ref{fig8} display  snapshots of 
the system with dimension $\d = 2$ and
birth rate 
\eqref{truncated_b_rate_k} with $k = 5$.

 {\textbf{5}}.
We also think that the speed of propagation
has superadditive structure. For a birth rate $b$
satisfying our assumptions,
let $s(b)$ be the speed of propagation.

\begin{conj}
 For any birth rates $b_1$, $b_2$ satisfying our assumptions,
 we have
 \[
  s(b_1) + s(b_2) \leq s(b_1 + b_2).
 \]

\end{conj}

\section{The construction and properties of the process}
\label{sec6}

Here we proceed to construct the process 
as a unique solution to a stochastic integral
equation.
First such a scheme was carried out by
Massouli{\'e} \cite{Mass98}.
This method can be deemed an 
analog of the construction from graphical representation.
We follow here \cite{Bez15}. 

\begin{rmk}
  Of course, the process starting from a fixed initial condition we consider here 
  can be constructed as the minimal jump process
  (pure jump type Markov processes 
  in the terminology of \cite{KallenbergFound})
  as is done for example in \cite{EiWag03}. 
  Note however that we use coupling of infinitely 
  many processes starting at different time points
  from different initial conditions, 
  so we here employ another method.
\end{rmk}

Recall that
 \[
 \Gamma _0(\R ^\d)=\{ \eta \subset \R ^\d : |\eta| < \infty \},
\]
and the $\sigma$-algebra on $\Go$ 
was introduced in \eqref{sigmaalgebra}.
To construct the family of processes $(\eta ^{q,A} _t)_{t \geq q}$,
we
consider the stochastic equation with Poisson noise
\begin{equation} \label{se}
\begin{aligned}
|\eta _t \cap B| = \int\limits _{(q,t] \times B \times [0, \infty )   }
& \mathds{1} _{ [0,b(x, \eta _{s-} )] } (u) 
N(ds,dx,du)
+ |\eta _q \cap B|,  \quad
t \geq q, \  B\in \mathscr{B}(\R ^\d),
\end{aligned}
\end{equation}
 where
$(\eta _t)_{t \geq q}$ is a cadlag $\Go$-valued solution
 process,
 $N$ is a  Poisson point process on 
$\R_+  \times \R^\d \times \R_+  $,
the mean measure of $N$ is $ds \times dx \times du$.
We require the
processes $N$ and $\eta _0$ to be independent of each other. 
Equation \eqref{se} is understood in the sense that the equality
holds a.s. for every bounded
$B \in \mathscr{B} (\R ^\d) $ and $t \geq q$.
In the integral on the right-hand side of \eqref{se},
$x$ is the location and 
$s$ 
is the time of birth of a new particle. Thus, 
the integral over $B$ from $q$ to $t$ represents the number 
of births inside $B$ which occurred before $t$.

Let us assume for convenience that $q = 0$.
We will make the following assumption on the initial  condition:
\begin{equation} \label{condition on eta _0}
 E |\eta _0| < \infty.
\end{equation}

  We say that the process
$N$ is  \textit{compatible} with 
an increasing, right-continuous 
and complete filtration of 
$\sigma$-algebras
 $(\mathscr{F}_t, t \geq 0)$
if $N$ is adapted, that is, all random variables
of the type $N(\bar T _1, U)$, $\bar T _1 \in \mathscr{B}([0;t])$,
$U \in \mathscr{B}(\R ^\d \times \R _+)$, are $\mathscr{F}_t$-measurable,
and all random variables of the type
$N(t+h,U) - N(t,U)$, $ h \geq 0 , U \in \mathscr{B}(\R ^\d \times \R _+)$,
are independent of $\mathscr{F}_t$, $N(t,U) = N([0;t], U)$.

 \begin{defi} \label{weak solution}
 A \emph{(weak) solution} of equation \eqref{se} is a triple 
 $(( \eta _t )_{t\geq 0} , N )$, $(\Omega , \mathscr{F} , P) $,
 $(\{ \mathscr {F} _t  \} _ {t\geq 0}) $, where 
\end{defi}
 
  (i) $(\Omega , \mathscr{F} , P)$ is a probability space, 
and $\{ \mathscr {F} _t  \} _ {t\geq 0}$ is an increasing, right-continuous
 and complete filtration of sub-$\sigma$-algebras of $\mathscr {F}$,

  (ii)   $N$ is a
 Poisson point process on $\R_+  \times \R^\d \times \R_+ $  with  
 intensity  $ds \times dx \times du $,

  (iii) $ \eta _0  $ is a random  $\mathscr {F} _0$-measurable
  element in $\Go$
  satisfying \eqref{condition on eta _0},

  (iv) the processes $N$ and $\eta _0$ are independent,
 $N$ is
  compatible with $\{ \mathscr {F} _t  \} _ {t\geq 0} $,

  (v) $( \eta _t )_{t\geq 0} $ is a cadlag $\Go$-valued process
adapted to $\{ \mathscr {F} _t  \} _ {t\geq 0} $, $\eta _t \big| _{t=0} = \eta _0$,  

  (vi) all integrals in \eqref{se} are well-defined, 
\[
E \int\limits _0 ^t ds \int\limits _{\R ^\d} b(x, \eta _{s-}) 
  < \infty, \ \ \ t > 0,
  \]
  
    (vii) equality \eqref{se} holds a.s. for all $t\in [0,\infty]$
  and all Borel sets $B$.

  Let 
 \begin{align}\label{nonchalant}
   \mathscr{S} ^{0} _t =  \sigma \bigl\{ & \eta_0 , 
 N([0,q] \times B \times C ) ,   \\ &
 q \in [0,t], 
B \in \mathscr{B} (\R^\d), C \in \mathscr{B} (\R_+)  \bigr\},
 \notag
 \end{align}
 and
 let $\mathscr{S} _t$ be the completion of $\mathscr{S} ^{0} _t$ under $P$.
 Note that $\{ \mathscr{S} _t \}_{t\geq 0} $ 
 is a right-continuous filtration 
 (see Remark \ref{nonchalant2}).

  \begin{defi} \label{def strong solution}
    A solution  of \eqref{se} is called \emph{strong}
 if $( \eta _t )_{t\geq 0} $ is adapted to 
$(\mathscr{S} _t, t\geq 0)$.
  \end{defi}

\begin{rmk}
 In the definition above we considered solutions as processes indexed by 
 $ t\ \in[0,\infty)$. 
 The reformulations for 
the case  $ t \in [0,T]$, $0<T<\infty$, are straightforward. 
This remark also applies to many of the results below.
\end{rmk}
  
\begin{defi} \label{joint uniqueness in law}
 We say that \textit{joint uniqueness in law} holds for equation \eqref{se} with an initial
distribution $\nu$ if any two (weak) solutions $((\eta_t) , N)$ and 
$((\eta_t  ^{ \prime }) , N  ^{\prime} )$ of \eqref{se},
$Law(\eta _0)= Law(  \eta _0  ^{\prime})=\nu$, have the same joint distribution:

$$Law ((\eta_t) , N)
= Law ((\eta_t  ^{\prime} ), N ^{\prime}  ).$$

\end{defi}
  
  \begin{thm}\label{core thm}
Pathwise uniqueness, strong existence and joint uniqueness
in law hold for equation \eqref{se}.
 The unique solution is a  Markov process. 
\end{thm}

  \textbf{Proof}. Without loss of generality assume that
 $P \{ \eta _0 \ne \varnothing \} = 1$.
 Define the sequence of random pairs $\{(\sigma _n, \zeta _{\sigma _n}) \}$, where
\[
\sigma _{n+1}= \inf \{ t>0 : \int\limits _{(\sigma _n,\sigma _n+t] \times B \times [0, \infty )  }
 \mathds{1}_{[0,   b(x,  \zeta _{\sigma _n})]} (u) N(ds,dx,du) >0 \}+ \sigma _n, \ \ \sigma _0 = 0,
\]
and
\[
 \zeta _{0} = \eta _0, \ \ \ 
\zeta _{\sigma _{n+1}} = \zeta _{\sigma _n} \cup \{ z_{n+1} \}
 \]
for $z_{n+1} = \{x\in \R^\d: N (\{\sigma _{n+1}\} \times \{x \} \times [0,   b(x,  \zeta _{\sigma _n})] 
) >0  \}$.
The points 
$z_n$ are uniquely determined a.s. 
Furthermore, 
 $\sigma _{n+1} > \sigma _n$ a.s., and $\sigma _n$
 are finite a.s by \eqref{harlot}.
We define $\zeta _t = \zeta _{\sigma _n}$
for $t\in [\sigma _n, \sigma _{n+1})$. Then by induction on $n$ 
it follows that $\sigma _n $ is a stopping time for each $n \in \N$, and 
$\zeta _{\sigma _n}$ is $\mathscr{F} _{\sigma _n} $-measurable. 
By direct substitution we see that
 $(\zeta _t)_{t\geq 0}$ is a strong
 solution to \eqref{se} on the time interval
$t\in [0, \lim\limits _{n\to \infty} \sigma _n)$. 
Although we have not defined what is a solution,
or a strong solution, on a 
random time interval, we do not discuss it here.
Instead we are going to show that 
\begin{equation}\label{staunch}
\lim\limits _{n\to \infty} \sigma _n = \infty \ \ \ \textrm{a.s.}
\end{equation}

The process 
$(\zeta _t)_{t\in [0, \lim\limits _{n\to \infty} \sigma _n)}$
has the Markov property, because 
the process $N$ has the strong Markov property and independent increments. 
Indeed, conditioning on $ \mathscr{I}_{\sigma_n}$,
\[
 E \bigl[   \mathds{1}_{\{\zeta _{\sigma _{n+1}} = \zeta _{\sigma _{n}} \cup x 
 \text{ for some } x \in B \}} \mid \mathscr{I}_{\sigma_n}\bigr]=
 \frac{\int\limits _{ B}   b(x, \zeta _{\sigma _{n}}) dx}{\int\limits _{\R^\d }
 {b}(x,\zeta _{\sigma _n}) dx },
\]
thus the chain $\{ \zeta _{\sigma _{n}}\}_{n \in Z_+}$ is a Markov chain,
and, given $\{ \zeta _{\sigma _{n}}\}_{n \in Z_+}$,
$\sigma _{n+1} - \sigma _n $ are distributed exponentially:
\[
E\{ \mathds{1}_{\{\sigma _{n+1} - \sigma _n >a\}}  \mid \{ \zeta _{\sigma _{n}}\}_{n \in Z_+}\} 
= \exp \{ - a \int\limits _{\R^\d }
 {b}(x,\zeta _{\sigma _n}) dx \}.
\]
Therefore, the random variables 
$\gamma _n = (\sigma _{n} - \sigma _{n-1}){\int\limits _{\R^\d }
 {b}(x,\zeta _{\sigma _n}) dx}$ 
constitute an independent 
of $\{ \zeta _{\sigma _{n}}\}_{n \in Z_+}$
sequence of independent unit exponentials.
 Theorem 12.18 in \cite{KallenbergFound} implies  that
$(\zeta _t)_{t\in [0, \lim\limits _{n\to \infty} \sigma _n)}$
is a pure jump type Markov process.

The jump rate of $(\zeta _t)_{t\in [0, \lim\limits _{n\to \infty} \sigma _n)}$
is given by
\[
 c(\alpha) = \int\limits _{\R^\d }
 {b}(x,\alpha) dx .
\]
Condition \ref{sublinear growth} implies that 
$c(\alpha) \leq   ||a||_{1} \cdot |\alpha|$,
where $||a||_{1} =||a||_{ L^1(\R ^\d)}$. Consequently,
\[
c (\zeta _{\sigma _n} ) \leq ||a||_{1} \cdot |\zeta _{\sigma _n}| = 
||a||_{1} \cdot |\eta _0| + n ||a||_{1} .
\]

We see that $\sum _n \frac{1}{c (\zeta _{\sigma _n} )} = \infty$ a.s., 
hence Proposition 12.19 in \cite{KallenbergFound} 
implies that $\sigma _n \to \infty$.

We have proved the existence of a strong solution. 
The uniqueness follows by induction on jumps of the process. 
Namely, let $( \tilde{\zeta} _t )_{t\geq 0} $ be 
a solution to \eqref{se}. From (vii)
 of Definition \ref{weak solution} and the equality

$$
\int\limits _{ (0,\sigma _1) \times \R ^\d  \times [0, \infty ] }
 \mathds{1} _{ [0, {b}(x, {\eta} _0 )] }(u) N(ds,dx,du) = 0 ,
$$
it follows
that $P \{ \tilde{\zeta} \text{\ has a birth before \ }
\sigma _1  \} = 0$. At the same time, the equality

$$
\int\limits _{ \{\sigma _1\} \times \R ^\d  \times [0, \infty ]  }
 \mathds{1} _{ [0, {b}(x, {\eta} _0 )] }(u) N(ds,dx,du) = 1 ,
$$
which holds a.s., 
yields that $\tilde{\zeta}$ too has a birth at the moment $\sigma _1$, and in the same point 
of space at that. Therefore, $\tilde{\zeta}$ coincides with $\zeta$ up to $\sigma _1$ a.s. Similar reasoning shows
that they coincide up to $\sigma _n$ a.s., and, since $\sigma _n \to \infty$ a.s., 

$$P \{ \tilde{\zeta} _t = {\zeta} _t \text{\ for all \ } t\geq 0 \} = 1.$$

 Thus, pathwise uniqueness holds. Joint uniqueness in law follows
 from the functional dependence between 
 the solution to the equation and the `input'
 $\eta _0$ and $N$.
\qed

\begin{prop} \label{turd}
 If $b$ is rotation invariant, then so is $(\eta _t)$.
\end{prop}

\textbf{Proof}. It is sufficient to note that 
$( M _ \d \eta _t)$, where $M _\d \in \textrm{SO}(\d)$, is the unique solution 
to \eqref{se} with $N$ replaced by $M _\d ^{-1} N$ defined by
\[
 M _\d ^{-1} N ([0,q] \times B \times C) = N ([0,q] \times M _\d ^{-1} B \times C), 
 \ \ \ q \geq 0, B \in \mathscr{B} (\R ^\d), C \in \mathscr{B} (\R _+).
\]
$M _\d ^{-1} N$ is a Poisson point process with the same intensity, therefore by
uniqueness in law $( M _ \d \eta _t) \overset{d}{=} (\eta _t)$.

\begin{prop}\label{strongMP} \textbf{(The strong Markov property)}
  Let $\tau$ be an $(\mathscr{S} _t, t\geq 0)$-stopping time
  and let $\tilde \eta _0  \overset{d}{=} \eta _{\tau}$. Then 
  \begin{equation} \label{flip out}
   (\eta _{\tau + t}, t \geq 0) \overset{d}{=}
   (\tilde \eta _{ t}, t \geq 0).
  \end{equation}
  Furthermore, 
  for any $\mathscr{D} \in \mathscr{B}(D _{\Go} [0, \infty))$,
  \begin{equation*}
  P\{ (\eta _{\tau + t}, t \geq 0) \in \mathscr{D}
   \mid \mathscr{S} _{\tau}\} = 
   P\{ (\eta _{\tau + t}, t \geq 0) \in \mathscr{D}
   \mid \eta _{\tau}\};
  \end{equation*}
that is,
given $\eta _{\tau}$,  $(\eta _{\tau + t}, t \geq 0)$
 is conditionally independent  of $(\mathscr{S} _t, t\geq 0)$.
\end{prop}

\textbf{Proof}. Note that

\begin{equation*}
\begin{aligned}
|\eta _{\tau + t} \cap B| = \int\limits _{(\tau, \tau + t] \times B \times [0, \infty )  }
& \mathds{1} _{ [0,b(x, \eta _{s-} )] } (u) 
N(ds,dx,du)
+ |\eta _ \tau \cap B| ,  \quad t \geq 0,
\  B\in \mathscr{B}(\R ^\d).
\end{aligned}
\end{equation*}

Since the unique solution is adapted
to the filtration generated by the noise 
and initial condition, the conditional independence follows,
and \eqref{flip out} follows from the uniqueness in law.
We rely here on the strong Markov property 
of the Poisson point process,
 see
 Proposition \ref{Strong MP} below.
\qed

\begin{cor}\label{cor69}
Let $\tau$ be an $(\mathscr{S} _t, t\geq 0)$-stopping time
and $\{y \}$ be an $\mathscr{S} _{\tau}$-
measurable finite random singleton.
Then 
\[
 (\eta ^{\tau, \{ y\}} _{\tau + t} - y)_{t\geq 0}
 \overset{(d)}{=} (\eta _{t})_{t\geq 0}.
\]
\end{cor}
\textbf{Proof}. This is a consequence of
Theorem \ref{core thm}
and Proposition \ref{strongMP}.
\qed

Consider two growth processes $(\zeta ^{(1)}) _t$
and $(\zeta ^{(2)}) _t$
defined on the common probability space
ans satisfying 
equations of the form \eqref{se},

\begin{equation} \label{2se}
\begin{split}
|\zeta ^{(k)} _t \cap B |= \int\limits _{(q,t] \times B \times [0, \infty ) }
& \lambda \mathds{1} _{ [0,b _k(x,\zeta ^{(k)} _{s-} )] } (u) 
N(ds,dx,du) 
+ |\zeta ^{(k)} _q \cap B|,
\ \  k=1,2.
\end{split}
\end{equation}

Assume that and the rates
$b_1$ and $b_2$ satisfy the conditions of 
imposed on $b$ in Section 2. Let
$(\zeta ^{(k)}_t)_{t \in [0,\infty )}$ 
be the unique strong solutions.  

\begin{lem} \label{couple}
 Assume that a.s. $\zeta ^{(1)}_0 \subset \zeta ^{(2)}_0$, and
for any two finite configurations $\eta ^1 \subset \eta ^2$,
\begin{equation} \label{culpable}
b_1(x,\eta ^1) \leq b_2 (x, \eta ^2), \ \ \ x\in \R^\d.
\end{equation}
Then a.s.
\begin{equation} \label{culprit}
 \zeta ^{(1)}_t \subset \zeta ^{(2)}_t, \ \ \ t \in [0,\infty ).
\end{equation}
\end{lem}

\textbf{Proof}. 
Let $(\sigma _n) _{n\in \N}$ be
the ordered sequence of
the moments of births for 
 $ (\zeta ^{(1)}_t)$,
that is,
$
 t \in (\sigma _n) _{n\in \N} 
$ if and only if 
$|\zeta ^{(1)}_t \setminus \zeta ^{(1)}_{t-}| = 1$.
It suffices to show that for each $n \in \N$, $\sigma_n$ is a moment of birth
for $(\zeta ^{(2)}_t)_{t \in [0,\infty )}$ too,
and the birth occurs at the same place.
We use induction on $n$.

 Here we deal only with
the base case, the induction step is done in the same way.
Assume that 
$$
\zeta ^{(1)}_{\sigma _1} \setminus \zeta ^{(1)}_{\sigma _1 - } =\{x_1 \}.
$$
The process $(\zeta ^{(1)})_{t \in [0,\infty )}$ satisfies
\eqref{2se}, therefore
$N(\{x \} 
\times [0,b _k(x_1,  \zeta ^{(1)} _{\sigma _1-} )]) = 1$.
Since 
$$ 
\zeta ^{(1)} _{\sigma _1-} = \zeta ^{(1)} _0
\subset \zeta ^{(2)} _0 \subset \zeta ^{(2)} _{\sigma _1-},
$$
by \eqref{culpable} 
$$
N_1(\{x \} \times \{ \sigma _1 \} \times [0,b _k(x_1,  \zeta ^{(2)} _{\sigma _1-} )]) = 1,
$$
hence 
$$
\zeta ^{(2)}_{\sigma _1} \setminus \zeta ^{(2)}_{\sigma _1 - } =\{x_1 \}.
$$
\qed

 \section{Appendix. The strong Markov property of a Poisson point process}

We need the strong Markov property of a Poisson point process.
Denote $X := \R ^\d \times \R _+$
(compare 
the proof of
Proposition \ref{strongMP}),
and let 
 $l$ be
  the Lebesgue measure on $X$.
Consider a 
 a Poisson point process $N$ on
$\R _+ \times X $
with
intensity measure $dt\times l$.
Let $N$ be compatible with a right-continuous complete filtration
$\{ \mathscr{F}_t \}_{t \geq 0}$, and $\tau$ be a finite a.s.
$\{ \mathscr{F}_t \}_{t \geq 0}$-stopping time . Introduce another
Point process $\overline N $ on $\R _+ \times X $,

\[
 \overline N ([0;s] \times U) = N ((\tau;\tau + s] \times U), 
 \ \ \  U \in \mathscr{B}(X)).
\]

\begin{prop}\label{Strong MP}
 The process $\overline N$ is a Poisson point process
 on $\R _+ \times X $
 with
 intensity $dt\times l$, independent of 
 $\mathscr{F}_{\tau}$.
\end{prop}

\textbf{Proof}. To prove the proposition,
it suffices to show that

(i) for any $b>a>0$ and open bounded $U \subset X$, 
$\overline N ((a;b),U)$ is a Poisson random variable
with mean $(b-a) l (U)$, and

(ii) for any $b_k>a_k>0$, $k=1,...,m$, and any open bounded
$U_k \subset X$, such that $((a_i;b_i) \times U_i) \cap( (a_j;b_j) \times U_j) = \varnothing $, $i \ne j$,
 the collection $\{\overline N ((a_k;b_k) \times U_k)\}_{k=1,m}$ 
 is a sequence of independent random variables, independent of $\mathscr{F}_{\tau}$.

Indeed, $\overline N $ is determined completely by values on sets of type
$(b-a) \beta (U)$, $a,b,U$ as in (i), therefore it must be
an independent of $\mathscr{F}_{\tau}$ Poisson
point process 
if (i) and (ii) hold.

Let $\tau _n$ be the sequence of $\{ \mathscr{F}_t \}_{t \geq 0}$-stopping times, 
$\tau _n = \frac{k}{2^n}$ on $\{ \tau \in (\frac{k-1}{2^n};\frac{k}{2^n}] \}$, 
$k \in \N$. Then
 $\tau _n \downarrow \tau$ and
$\tau _n - \tau \leq \frac{1}{2^n}$. Note that the
stopping times $\tau _n$ take countably many values only.
The process $N$ satisfies the strong Markov property for $\tau _n$: 
 the processes $\overline N _n$, defined by
$$
\overline N _n ([0;s] \times U) := N ((\tau_n;\tau_n + s] \times U), 
$$
are Poisson point processes, independent of $\mathscr{F}_{\tau _n}$.
To prove this, take $k$ with $P\{ \tau_n = \frac{k}{2^n} \} >0$ and
note that on $\{ \tau_n = \frac{k}{2^n} \}$,  $\overline N _n $ 
coincides with
 process the Poisson point process $\tilde N _{\frac{k}{2^n}}$
given by
\[
\tilde N _{\frac{k}{2^n}} ([0;  s] \times U)
:=  N  \bigg((\frac{k}{2^n};\frac{k}{2^n} + s] \times U) \bigg), \ \ \  U \in \mathscr{B}(\R ^\d).
\]
Conditionally on $\{ \tau_n = \frac{k}{2^n} \}$,
$\tilde N _{\frac{k}{2^n}}$ is again a Poisson point process,
with the same intensity. Furthermore, conditionally on $\{ \tau_n = \frac{k}{2^n} \}$,
$\tilde N _{\frac{k}{2^n}}$
is independent of $\mathscr{F}_{\frac{k}{2^n}}$, hence it
is independent 
of 
$ \mathscr{F}_{\tau } \subset \mathscr{F}_{\frac{k}{2^n}}$.

To prove (i), 
note that $\overline N _n ((a;b) \times U) \to \overline N ((a;b) \times U) $ a.s. 
and all random variables $\overline N _n ((a;b) \times U)$ have the same distribution, 
therefore $\overline N  ((a;b) \times U)$ is a Poisson random variable with mean 
$(b-a)\lambda (U)$. The random variables $\overline N _n ((a;b) \times U)$
are independent of $\mathscr{F}_{\tau}$, 
hence $\overline N  ((a;b) \times U)$ is independent 
of $\mathscr{F}_{\tau}$, too. Similarly, (ii) follows. $\Box$

\begin{rmk}\label{nonchalant2}
We assumed in   
Proposition \ref{Strong MP}
that there exists an increasing, 
 right-continuous and complete 
filtration 
$\{ \mathscr{S}_t \}_{t \geq 0}$
compatible with $N$. Let us show
 that such  filtrations exist.

Introduce the natural filtration of $N$,

\[
 \bar {\mathscr{S}}_t ^0 = \sigma \{  N_k(C, B), 
B\in  \mathscr{B} (\R ^\d), C\in  \mathscr{B} ([0;t]) \},
\]
and let $\bar {\mathscr{S}}_t $ be 
the completion of $\bar {\mathscr{S}}_t ^0$
under $P$.
Then $N$ is compatible with $\{\bar {\mathscr{S}}_t\} $.
We claim that $\{\bar {\mathscr{S}}_t \}_{t \geq 0}$, 
defined in such a way,
is right-continuous (this may be regarded as an 
analog of Blumenthal's 0-1 law). Indeed, as in the proof 
of Proposition \ref{Strong MP}, we can  check
that $\tilde N _a$ is independent of $\bar {\mathscr{S}}_{a+}$.
Since 
$\bar {\mathscr{S}}_{\infty}  = \sigma(\tilde N _a) \vee 
\bar {\mathscr{S}}_{a}$,
$\sigma(\tilde N _a)$ and $\bar{\mathscr{S}}_{a}$ are 
independent
and $\bar {\mathscr{S}}_{a+}
\subset  \bar {\mathscr{S}}_{\infty} $,
we see that
$\bar{\mathscr{S}}_{a+} \subset\bar {\mathscr{S}}_{a} $. Thus, 
$\bar{\mathscr{S}}_{a+} =\bar {\mathscr{S}}_{a} $.

\end{rmk}

 \section*{Acknowledgement}
 
 Luca Di Persio would like to acknowledge the Informatics Department at
 the University of Verona for having funded the project
``Stochastic differential equations with jumps in 
Mathematical Finance: applications to pricing, hedging and dynamic risk measure's problems''
 and the Gruppo Nazionale per l'analisi matematica, la Probabilit\`{a} e le loro applicazioni (GNAMPA).
  Viktor Bezborodov is  supported by
the Department of Computer Science at the University of Verona.
Viktor Bezborodov was also partially supported by the DFG through the SFB 701
``Spektrale Strukturen und Topologische Methoden in der Mathematik'' and
by the European Commission under the project STREVCOMS PIRSES-2013-612669.
Mykola Lebid is supported by the Department of Biosystems Science and Engineering
at the ETH Z\"urich. Tyll Krueger
and Tomasz O\.za\'nski are supported by the Wroclaw University of Technology.

\bibliographystyle{alpha}
\bibliography{Sinus}

\end{document}